\documentclass[12pt]{amsart}
\usepackage{amsmath,amssymb,latexsym,esint,cite,mathrsfs}
\usepackage{verbatim,wasysym,cite}
\usepackage[left=2.8cm,right=2.8cm,top=2.9cm,bottom=2.9cm]{geometry}
\usepackage{xcolor}

\usepackage[utf8]{inputenc}

\usepackage{tikz}
\usepackage{pgfplots}
\pgfplotsset{compat=1.11}

\usepackage[colorlinks=true,urlcolor=blue, citecolor=red,linkcolor=blue,
linktocpage,pdfpagelabels, bookmarksnumbered,bookmarksopen]{hyperref}

\numberwithin{equation}{section}
\newtheorem{theorem}{Theorem}[section]
\newtheorem{lemma}[theorem]{Lemma}
\newtheorem{remark}[theorem]{Remark}
\newtheorem{definition}[theorem]{Definition}
\newtheorem{proposition}[theorem]{Proposition}
\newtheorem{corollary}[theorem]{Corollary}
\DeclareMathOperator{\eq}{eq}

\DeclareMathOperator{\ran}{ran}
\DeclareMathOperator{\dom}{dom}
\DeclareMathOperator{\re}{Re}
\DeclareMathOperator{\im}{Im}
\DeclareMathOperator{\sign}{sign}

\DeclareMathOperator{\ess}{ess} 
\DeclareMathOperator{\dist}{dist} 
\DeclareMathOperator{\Span}{span}
\newcommand{\mc}[1]{\mathcal{#1}}

\title[ NLS-$\delta'_s$ equation on the star graph  ]
{Dynamical and variational properties of the NLS-$\delta'_s$ equation on the star graph}

\author[Nataliia Golshchapova]{}
\email{ nataliia@ime.usp.br }


\begin{document}

\maketitle

\centerline{Nataliia Goloshchapova}
{\footnotesize
\centerline{Institute of Mathematics and Statistics, University of S\~ao  Paulo,}
\centerline{Rua do Mat\~ao, 1010,  S\~ao  Paulo-SP, 05508-090, Brazil}
}

\begin{abstract}
We study the nonlinear Schr\"odinger equation with  $\delta'_s$ coupling of intensity $\beta\in\mathbb{R}\setminus\{0\}$ on the star graph $\Gamma$ consisting of $N$ half-lines. The nonlinearity   has the form $g(u)=|u|^{p-1}u, p>1.$ In the first part of the paper, under certain restriction on $\beta$, we prove the existence of the ground state solution as a minimizer of the action functional $S_\omega$ on the  Nehari manifold. It appears that the family of critical points which contains a ground state consists of $N$ profiles   (one symmetric and $N-1$ asymmetric). In  particular, for the attractive $\delta'_s$ coupling ($\beta<0$) and  the frequency $\omega$ above a certain threshold,  we managed to specify the ground state. 

  The second part is devoted to the study  of  orbital instability   of the  critical points.  We  prove spectral instability of the critical points using Grillakis/Jones Instability Theorem. Then  orbital instability for $p>2$  follows from the fact that data-solution mapping associated with the equation  is of class $C^2$ in $H^1(\Gamma)$. Moreover, for $p>5$ we complete and concertize instability results showing  strong instability (by blow up in finite time) for the particular  critical points.

\end{abstract}
\section{Introduction}
\label{}
 We consider a star graph $\Gamma$ consisting of a central vertex $\mathbf{v}$ and $N\geq 2$ infinite half-lines attached to it.  We identify $\Gamma$ with the disjoint union of the intervals $I_j=(0,\infty)$, $j=1,\ldots,N$, augmented by the central vertex $\mathbf{v}=0$.  The function on $\Gamma$ is defined by
  $$v:\Gamma\rightarrow\mathbb{C}^{N},\quad v=(v_j)^{N}_{j=1},\quad v_j:(0,\infty)\rightarrow\mathbb{C}.$$
   We  denote by $v_j(0)$ and $v'_j(0)$ the limits of $v_j(x)$ and  $v'_j(x)$  as $x\rightarrow 0+$. 

The principal object of this study is the following focusing nonlinear Schr\"odinger equation on  $\Gamma$ with $\delta'_s$ coupling:
\begin{equation}\label{NLS}
\begin{cases} 
i\partial_{t}u(t,x)=-\Delta_{\beta}u(t,x)-\left|u(t,x)\right|^{p-1}u(t,x),\quad(t,x)\in \mathbb{R}\times\Gamma,\\
u(0,x)=u_{0}(x),
\end{cases} 
\end{equation}
where $\beta\in \mathbb{R}\setminus\{0\}$, $p>1$, $u:\mathbb{R}\times\Gamma\rightarrow\mathbb{C}^{N}$,  and  $\left(-\Delta_{\beta}v\right)(x)=(-v_j''(x)), x>0$, with  
\begin{equation}\label{dom_delta'}
\dom(\Delta_{\beta})=\Big\{v\in H^1(\Gamma): v'_1(0)=\ldots =v'_N(0),\quad \sum^{N}_{j=1}v_j(0)=\beta v'_{1}(0)\Big\}.
\end{equation} 
The Schr\"odinger operator  $ -\Delta_\beta$ with $\delta'_s$ coupling has the precise interpretation as the self-adjoint operator on $L^{2}(\Gamma)$ uniquely associated (by the KLMN Theorem \cite[Theorem X.17]{ReeSim78}) with the closed semibounded quadratic form $ {F}_{\beta}$ defined on $H^{1}(\Gamma)$ by 
\begin{equation*}\label{F}
 {F}_{\beta}(v)=\|v'\|^{2}_2+\dfrac{1}{\beta}\left|\sum\limits_{j=1}^N v_j(0)\right|^{2}.
\end{equation*}
Similarly to the case of $\delta'$ coupling on the line (see \cite{AvrExn94}),  the $\delta'_s$ coupling on $\Gamma$ has a high energy scattering behavior that can be reproduced through scaling limits of scatterers, the so-called spiked-onion graphs. These are obtained replacing $\Gamma$ by 
the graph $\Gamma_{N}(n, \ell)$, where  every pair of half-line endpoints is connected by $n$ links of length $\ell$ (see \cite{Exn96}). The corresponding variables  run through the interval $[-\ell / 2, \ell / 2] ,$  and it is assumed that  the coupling at each graph vertex is $\delta$ (see formula \eqref{delta}).
One assumes that  $\ell \rightarrow 0, n \rightarrow \infty$ while keeping the product $n \ell$ fixed.

In \cite{CheExn04} the authors proposed an interpretation of the Hamiltonian $-\Delta_\beta$  as a norm resolvent limit of the operator $H(\varepsilon)$ when  $\varepsilon\to 0+$. 
The operator $H(\varepsilon)$   acts as  minus second derivative on each edge of   the graph $\Gamma_{\varepsilon}$ consisting of $\Gamma$ with additional vertices of degree two at each half-line,  all at the same distance $\varepsilon>0$ from  $\mathbf{v}=0$.  The wave functions satisfy the $\delta$ coupling conditions  at   $\mathbf{v}$ and at new vertices of degree two.  The intensity of  $\delta$ couplings depends on $\beta$ and $\varepsilon$.

The operator $-\Delta_\beta$ belongs to $N^2$-parametric family of self-adjoint extensions of  the minimal symmetric operator 
\begin{align*}
&\left(-\Delta_{\mathrm{min}} v\right)(x)=\left(-v_{j}^{\prime \prime}(x)\right)_{j=1}^{N}, \quad x>0, \quad v=(v_{j})_{j=1}^{N},\\
&\dom\left(\Delta_{\mathrm{min}}\right)=\left\{v \in H^{2}(\Gamma): v_{1}(0)=\ldots=v_{N}(0)=0,\quad v'_{1}(0)=\ldots=v'_N(0)=0\right\}.
\end{align*}
The whole  family of self-adjoint conditions naturally arising at the vertex $\mathbf{v}=0$ of  $\Gamma$ has the  description (see \cite{CheExn04}):
$$
(U-I)v(0) +i(U+I)v'(0)=0,
$$
where $v(0)=\left(v_{j}(0)\right)_{j=1}^{N}, v'(0)=\left(v'_j(0)\right)_{j=1}^{N}, U$ is an arbitrary
unitary $N \times N$ matrix, and $I$ is the $N \times N$ identity matrix.
In our case $U=I-\frac{2}{N-i \beta} \mathcal{I}$, where $\mathcal{I}$ is  $N\times N$  matrix whose all entries equal
one.
It is a difficult problem to understand which of self-adjoint conditions are physically relevant (self-adjointness is just a necessary physical requirement to ensure conservation of the probability
current at the vertex).

Among all the possible matching conditions, the
most used are the  Kirchhoff ones:
$$v _1(0)=\ldots =v_N(0),\quad \sum^{N}_{j=1}v'_j(0)=0.$$
Justifications
of the Kirchhoff conditions on different types of metric graphs have been obtained
in many physical experiments involving wave propagation in thin waveguides
and quantum nanowires. Namely, they appear when  multi-dimensional models are  approximated by differential operators 
  on graphs. Moreover, when imposing the Kirchhoff conditions,
transmission and reflection coefficients are independent of the momentum.

These arguments led to the fact that the Kirchhoff conditions have been assumed as the most natural, and 
hence they become the most widely studied. However, it is not clear whether these conditions  suit  different physical models (especially in the presence of
non-trivial localized interactions near junctions).  For instance, by choosing  different values of the thickness parameters vanishing
at the same rate, it was shown in \cite{Kuc01} that generalized Kirchhoff boundary conditions (or ``weighted'' Kirchhoff conditions) can also arise in the asymptotic limit.
Moreover, in \cite{CheTsu02, FulTsu00} 
the authors suggested that  
 some other conditions could be more
satisfactory from the point of view of invariance laws.
However, it is worth mentioning  that $\delta'_s$ coupling  cannot be achieved in a purely geometrical way, by squeezing a system of branching tubes with the same topology as the graph (see \cite{CheExn04} and  references therein). Moreover, it still  cannot be obtained using approximations involving potentials scaled in the usual way. Nevertheless,  $\delta'_s$ coupling being modeled  by complicated enough geometric scatterers, is likely to have something in common with the real world when one replaces simple junctions  by regions
of a nontrivial topology (see \cite[Section 4]{Exn96}).

The systematic  study of nonlinear  evolution equations on metric graphs dates back to \cite{Mah94}.
 The nonlinear
PDEs on graphs, mostly the nonlinear Schr\"odinger equation (NLS), have been studied in
the past decade in the context of existence, stability, and propagation of solitary waves (see \cite{Noj14} for the references).
Two  fields where NLS equation appears as a preferred model are optics of nonlinear
Kerr media and dynamics of Bose-Einstein condensates involving application to graph-like structures. For example, in nonlinear optics one studies 
arrays of planar self-focusing waveguides and propagation in variously shaped fiber optics devices,
such as $Y$-junctions, $H$-junctions  (see, for instance, \cite{BulSad11, HunTri11}).

The extensive study of  existence of ground states (as minimizers of energy functional under fixed mass) for the NLS models on metric graphs was carried out in the presence of the Kitchhoff conditions at the vertices of the graph (see \cite{AdaSer15, AdaSer16} and references therein). The authors analyzed the problem for metric graphs with
a finite number of vertices and at least one half-line. The existence and the stability of  solitary waves for different types of graphs with the Kirchhoff conditions at the vertices were treated in numerous papers \cite{KaiPel18a, KaiPel21, MarPel, PelSch17}. It is worth mentioning recent paper  \cite{NojPel20}, where  the authors  explored the variational methods and the analytical
theory for differential equations in order to construct the standing waves for the quintic NLS with Kirchhoff conditions on the tadpole graph.

Rigorous study of the NLS models on graphs in presence of impurities is related to a so-called $\delta$ coupling. On $\Gamma$ it is defined by:
\begin{equation}\label{delta}
v _1(0)=\ldots =v_N(0),\quad \sum^{N}_{j=1}v'_j(0)=\alpha v_1(0),\,\,\alpha\in\mathbb{R}\setminus\{0\}.
\end{equation}
  The $\delta$ coupling is the most studied non-Kirchhoff condition \cite{AQFF, AdaNoj14,  AnGo2018, BanIgn14, CAFiNo2017,  NataMasa2020, Gol21, Kai19}. 
In \cite{AQFF, AdaNoj14}, for $\alpha<0$, the ground state solution (as a minimizer of action functional and energy functional, respectively) has been identified with the $N$-tail state only in presence of sufficiently  strong attractive interaction and sufficiently small mass, respectively. 
Besides the $N$-tail symmetric state, the NLS with $\delta$ coupling on $\Gamma$  admits the family of asymmetric states which are constituted by  tail- and bump-components on the edges of $\Gamma$. Extensive study of their orbital stability was made in  \cite{AnGo2018, AngGol18a, Kai19}.

If we decide to drop the continuity condition, the next more general class of self-adjoint conditions which seems natural for applications consists of those which are locally permutation invariant at each vertex. 
The $\delta'_s$ coupling belongs to this class. Observe that it is a natural generalization of the symmetrized
version of  $\delta'$ coupling on the line (see \cite{AlbGes05} for the  comprehensive treatment of this coupling). 
Systematic  investigation of the  NLS model with the $\delta'$ coupling on the line appears in \cite{AdaNoj13}.
The authors prove the existence of the minimizer of the action functional $S_\omega$ on the Nehari manifold for   attractive coupling. 
It appears that the $\delta'$ coupling gives rise to a much
richer structure of the family of ground states, including a pitchfork bifurcation with
symmetry breaking. 
More precisely, there exists a critical value $\omega^{*}$ of frequency  such that if $\omega<\omega^{*}$, then there is a single ground state and it is an odd function; if $\omega>\omega^*$, then there exist two non-symmetric ground states (contrarily
to what happens in the case of the  $\delta$ coupling on the line, where all ground states are even functions).
This complex picture  of the ground states originates  from the fact that an associated 
energy space  includes functions discontinuous at zero. The investigation of orbital stability of the odd ground state in the case of repulsive $\delta'$ coupling was carried out in \cite{AngGol20}.

 In this paper we generalize  and extend the results obtained in \cite{AdaNoj13, AngGol20} for the NLS model on $\Gamma$. The study of NLS-$\delta'_s$ equation on $\Gamma$ was initiated in \cite{AnGo2018}. In the case of attractive coupling, the authors investigated orbital stability for  standing wave $e^{i\omega t}\phi(x)$ with the  $N$-tail profile $\phi$ (see \cite[Theorem 1.2]{AnGo2018}). Mathematically, the main
advantage of studying  $\delta'_s$ coupling   is the existence of an explicit nontrivial family of  soliton profiles to be described below. 
 
The present  manuscript has four principal parts. Firstly, we give detailed proof on well-posedness of   \eqref{NLS} in $H^1(\Gamma)$ (which is the energy space). In particular, in Proposition \ref{PP1} we show that for $1<p<5$ global well-posedness holds, and  for $p>2$  data-solution mapping associated with the equation  is of class $C^2$ in $H^1(\Gamma)$ (which is crucial for our proof of orbital instability). Additionally, we have shown the regularity result for $p>2$  and $u_0\in\dom(\Delta_{\beta})$ (see Proposition \ref{well_D_H}). This regularity result is important for the proof of virial identity \eqref{well_17a}.

Secondly, we  deal with the existence  of the ground states as minimizers of the action functional $S_\omega$ restricted to the Nehari manifold (see \eqref{d}). In Theorem \ref{GS} we prove that for attractive and sufficiently weak $\delta'_s$ coupling the minimizer exists. The principal  step in the proof is to compare  our minimization  problem with the one for $\beta=\infty$. This problem involves the technique  of   symmetric rearrangements on $\Gamma$ elaborated in \cite{AQFF}: that is, we are able to reduce the problem to the  half-line.
It is worth mentioning that in the space of symmetric functions $H^1_{\eq}(\Gamma)$ the minimizer exists for any $\beta\in\mathbb{R}\setminus\{0\}$ (Lemma \ref{lemma_eq}). Furthermore, for even  $N$ and $\beta<0$,  the minimizer always exists in  $H^1_{\frac{N}{2}}(\Gamma)$ as well (Lemma \ref{N/2}).  This situation is analogous to the case of the real line considered in \cite{AdaNoj13} (indeed, in this case $\Gamma$ can be identified with $\frac{N}{2}$ copies of $\mathbb{R}$).  In particular, for $\omega\leq\frac{p+1}{p-1}\frac{N^2}{\beta^2}$ the minimizer is given by the symmetric profile $\phi_\beta$, and for $\omega>\frac{p+1}{p-1}\frac{N^2}{\beta^2}$ it coincides with the  asymmetric  profile $\phi_{\frac{N}{2}}$ (with equal first $\frac{N}{2}$ entries and last $\frac{N}{2}$ entries).

 Thirdly, we are looking for the candidates to be the minimizers, i.e. critical    points to $S_\omega$ of the form $e^{i\theta}\phi(x)$, where $\phi(x)$ is a real-valued  profile. It appears that for $\omega>\frac{p+1}{p-1}\frac{N^2}{\beta^2}$ the whole family of such critical points consists of $N$ profiles modulo permutations of the edges of $\Gamma$ (one is symmetric $\phi_\beta$ and $N-1$ are asymmetric profiles $\phi_k$). We conjecture that for $\omega$ below $\frac{p+1}{p-1}\frac{N^2}{\beta^2}$ the symmetric profile $\phi_\beta$ is the minimizer, while in Theorem \ref{min_phi_1} we managed to prove that for $\omega>\frac{p+1}{p-1}\frac{N^2}{\beta^2}$ the minimizer is given by the  asymmetric  tail-profile $\phi_1$.
Notice that in \cite{AdaNoj13} the authors proved minimizing property of the asymmetric  profile by direct comparison of the values of the action functional (at symmetric and asymmetric profile). Our proof is completely different. We use the Implicit Function Theorem  and the fact that the Morse index of $S''_\omega(\phi_1)$ equals 1.
 
Lastly, we study instability properties of the family of the critical points mentioned above. Namely, using the  Grillakis/Jones Instability   Theorem  (see \cite{Grill88, Jon88}),
we have proved spectral instability (see Theorem \ref{main_asym})  of the asymmetric  critical points $\phi_k$ for $\beta<0, k\geq 2$ and $\beta>0, N-k\geq 4$.   
For $p>2$, using $C^2$ regularity  of the  mapping data-solution and applying the abstract result from \cite{HenPer82}, we have shown orbital instability of $\phi_k$. This abstract result  states  the nonlinear instability of a fixed point of a nonlinear mapping having the linearization $L$ of spectral radius $r(L)>1$.
To apply the  Grillakis/Jones approach we need to estimate the Morse index of two self-adjoint in $L^2(\Gamma)$ operators associated with $S''_\omega(\phi_k)$. These estimates were obtained in Proposition \ref{Morse_index_assym}  by using the generalization of the Sturm theory elaborated for $\Gamma$ in  \cite{KaiPel18, Kai19}.

Applying the  Grillakis/Jones Instability   Theorem  to the  symmetric profile $\phi_\beta$, in Theorem \ref{main_instub} we obtain  partial generalization of instability results \cite[Theorem 6.13]{AdaNoj13} and \cite[Theorem 1.1]{AngGol20}. 
Finally,  we concertize the instability  results by  claiming  strong instability (by blow up)  of  the symmetric profile $\phi_\beta$  in the supercritical case $p>5$ (see Theorem \ref{beta>0} and \ref{beta<0}). The proof essentially uses variational characterization of $\phi_\beta$ in $H^1_{\eq}(\Gamma)$ and virial identity  \eqref{well_17a}.

The paper has the following structure. In Section 2  we prove  well-posedness of problem \eqref{NLS}  in $H^1(\Gamma)$ and in $\dom(\Delta_\beta).$ Section 3 is devoted to the proof of the existence of a ground state solution. In Section 4 we find the family of critical points that contains  ground states. In Section 5 we provide an extensive study of the instability of standing  wave solutions associated with the mentioned critical points. Namely, Subsection \ref{symm_prof} is devoted to the case of symmetric profile $\phi_\beta$, while in Subsection \ref{assym} we study instability of the asymmetric profiles $\phi_k.$ 
 
\vspace{1cm}

\textbf{Notation.}
The natural Hilbert space associated to  the Laplace operator $-\Delta_{\beta}$  is $L^{2}(\Gamma)$, which is
defined as \\ $L^{2}(\Gamma)=\bigoplus^{N}_{j=1}L^{2}(\mathbb{R}^{+})$. The inner product in $L^{2}(\Gamma)$ is given by $$(u,v)_2=\re\sum\limits_{j=1}^N(u_j, v_j)_{L^2(\mathbb{R}^+)}, \quad u=(u_j)_{j=1}^N,\,v=(v_j)_{j=1}^N .$$
The space $L^{q}(\Gamma)$, for $1\leq q \leq\infty, $ is defined as prime sum as well, and  $\|\cdot\|_q$ stands for its norm.
By  $H^{1}(\Gamma)=\bigoplus^{N}_{j=1}H^{1}(\mathbb{R}^{+})$ and $H^{2}(\Gamma)=\bigoplus^{N}_{j=1}H^{2}(\mathbb{R}^{+})$ we define the  Sobolev spaces.
In what follows we  will use the notation  $H^1$ for $H^1(\Gamma)$.
The duality pairing between $(H^1(\Gamma))'$ and $H^1(\Gamma)$ is denoted by $\langle\cdot, \cdot\rangle_{(H^1)'\times H^1}.$

 Below we will run into the following subspaces of $H^1(\Gamma)$:
 \begin{align*}
&H_{\mathrm{eq}}^{1}(\Gamma)=\left\{v \in H^{1}(\Gamma): v_{1}(x)=\ldots=v_{N}(x), x>0\right\},\\
&H_{\frac{N}{2}}^1(\Gamma)= \left\{v \in H^{1}(\Gamma): v_{1}(x)=\ldots=v_{\frac{N}{2}}(x),\,\, v_{\frac{N}{2}+1}(x)=\ldots=v_{{N}}(x), x>0\right\}.
 \end{align*}
 
By $C_{j}, C_{j}(\cdot), K_j(\cdot), j \in \mathbb{N}$, and $C(\cdot)$ we will denote some positive constants.
Let $L$ be a self-adjoint operator in $L^2(\Gamma)$, then $n(L)$  will stand for the number of negative eigenvalues of $L$ counting multiplicity (Morse index). 
\iffalse
Note that we can formally rewrite  \eqref{NLS} as
\begin{equation*}
i\partial_{t}u(t)=E^{\prime}(u(t)),
\end{equation*}
where $E$ is the energy functional defined by
\begin{align*}
E(u)=\frac{1}{2} {F}_{\beta}(u)-\frac{1}{p+1}\|u\|^{p+1}_{p+1}.
\end{align*}
\fi
\section{Well-posedness}
In this section we study well-posedness of problem \eqref{NLS}.
The proposition below states  local well-posedness of \eqref{NLS} in $H^1(\Gamma).$
 \begin{proposition}\label{PP1}
Let $u_{0}\in H^{1}(\Gamma)$ and $p>1$. Then the following assertions hold. 

$(i)$ There exist $T=T(u_{0})>0$ and a unique weak solution $u(t)\in C\left([0,T], H^{1}(\Gamma)\right)\cap C^1\left([0,T], (H^{1}(\Gamma))'\right)$ of problem \eqref{NLS}.

$(ii)$ Problem \eqref{NLS} has a maximal solution defined on
an interval  $[0, T_{H^1})$, and the following ``blow-up alternative'' holds: either $T_{H^1} = \infty$ or $T_{H^1} < \infty$ and
$$\lim\limits_{t\to T_{H^1}}\|u(t)\|_{H^1} =\infty.$$

$(iii)$
For each $T_0\in (0, T)$ the mapping
$u_0\in H^1(\Gamma)\mapsto u(t)\in C\left([0, T_0], H^1(\Gamma)\right)$ is continuous. 
In particular, for $p>2$ this mapping  is at least of class $C^2$.

$(iv)$ The conservation of energy and charge holds:  for $t\in[0, T_{H^1})$
\begin{align}\label{conserv}
E(u(t))=\frac{1}{2} {F}_{\beta}(u(t))-\frac{1}{p+1}\|u(t)\|^{p+1}_{p+1}=E(u_{0}),\quad\left\|u(t)\right\|^2_2=\left\|u_{0}\right\|^2_2.
\end{align}
\end{proposition}
\begin{proof}
The proof is based on the ideas of the one of  \cite[Theorem 4.10.1]{Caz03}. 

\noindent\textit{Step 1.} To start, we mention several useful technical facts. Firstly, observe that  $g(u)=|u|^{p-1}u\in C^1(\mathbb{C},\mathbb{C})$ (i.e. $\im(g)$ and $\re(g)$  are $C^1$-functions of  $\re u, \im u$) for $p>1$.
This implies inequalities (see \cite[Lemma 4.10.2]{Caz03})
\begin{equation}\label{g_1}
\begin{aligned}
\|g(u)\|_{H^1} & \leq C_1(M)\|u\|_{H^1}, \\
\|g(u)-g(v)\|_{2} & \leq C_2(M)\|u-v\|_{2},
\end{aligned}
\end{equation}
for all $u, v \in H^1(\Gamma)$ such that $\|u\|_{\infty},\|v\|_{\infty} \leq M,$ and 
\begin{equation}\label{g_2}
\|g(u)-g(v)\|_{H^1} \leq C_3(M)\Big[\|u-v\|_{H^1}+\varepsilon_{M}\left(\|u-v\|_{2}\right)\Big]
\end{equation}
for all $u, v \in H^1(\Gamma)$ such that $\|u\|_{H^1},\|v\|_{H^1} \leq M,$ where $\varepsilon_{M}(s) \rightarrow 0$ as $s \downarrow 0$.

 Secondly, we show the   inequality 
\begin{equation}\label{group}
  \|e^{i \Delta_\beta t}v\|_{H^1}\leq C \|v\|_{H^1}.  
\end{equation} 
 
Let $m>\frac{N^2}{\beta}$. It is known that $\inf \sigma(-\Delta_\beta)=\left\{\begin{array}{c}0,\,\, \beta \geq 0, \\ -\frac{N^{2}}{\beta^{2}}, \beta<0 .\end{array}\right.$ Then  the operator $-\Delta_\beta+m$ is positive. 
Remark that $H^1(\Gamma)=\dom\left( {F}_{\beta}\right)=\dom((-\Delta_\beta+m)^{1/2})$.  Using $L^2$-unitarity of $e^{i\Delta_\beta t}$, we obtain for $v\in H^1(\Gamma)$ 
\begin{align*}
 & {F}_{\beta}(v)+m\|v\|_2^2=\left(( -\Delta_\beta+m)^{1/2}v,( -\Delta_\beta+m)^{1/2}v\right)_2\\&=\left(e^{i \Delta_\beta t}( -\Delta_\beta+m)^{1/2}v, e^{i \Delta_\beta t}( -\Delta_\beta+m)^{1/2}v\right)_2\\&=\left(( -\Delta_\beta+m)^{1/2}e^{i\Delta_\beta t}v, ( -\Delta_\beta+m)^{1/2}e^{i\Delta_\beta t}v\right)_2= {F}_\beta(e^{i\Delta_\beta t}v)+m\|e^{i\Delta_\beta t}v\|_2^2. 
\end{align*}
Using equivalence of $H^1$-norm and ${F}_{\beta}(v)+m\|v\|_2^2$ (see \cite[Lemma 4.13]{ArdCel21}),  we get
\begin{align*}
 &C_2\|e^{i\Delta_\beta t}v\|^2_{H^1}\leq  {F}_{\beta}(e^{i\Delta_\beta t}v)+m\|e^{i\Delta_\beta t}v\|_2^2 = {F}_{\beta}(v)+m\|v\|_2^2\leq C_1\|v\|^2_{H^1}, 
\end{align*}
and \eqref{group} follows easily.  
 Analogously one proves that $e^{i\Delta_\beta t}$ is continuous in $H^1(\Gamma).$
  Indeed, let $t_n\underset{n\to\infty}{\longrightarrow} t$, then  
  \begin{align*}
 &C_2\|(e^{i\Delta_\beta t}-e^{i\Delta_\beta t_n})v\|^2_{H^1}\leq  {F}_{\beta}((e^{i\Delta_\beta t}-e^{i\Delta_\beta t_n})v)+m\|(e^{i\Delta_\beta t}-e^{i\Delta_\beta t_n})v\|_2^2\\& = \Bigl((e^{i\Delta_\beta t}-e^{i\Delta_\beta t_n})( -\Delta_\beta+m)^{1/2}v,(e^{i\Delta_\beta t}-e^{i\Delta_\beta t_n})( -\Delta_\beta+m)^{1/2}v\Bigr)_2\\
 &=2\|(-\Delta_\beta+m)^{1/2}v\|^2_2-\Bigl((e^{-i\Delta_\beta(t-t_n)}+e^{i\Delta_\beta(t_n-t)})( -\Delta_\beta+m)^{1/2}v, ( -\Delta_\beta+m)^{1/2}v \Bigr)_2\\
 &\underset{n\to\infty}{\longrightarrow} 0
\end{align*}
\textit{Step 2.} Existence of  solution  follows  by a fixed point argument. Given $M, T>0$ to be chosen later, we set 
$$
E_1=\left\{u \in L^{\infty}\left((0,T), H^1(\Gamma)\right):\|u\|_{L^{\infty}\left((0,T), H^{1}(\Gamma)\right)} \leq M\right\},
\quad
{d}_1(u, v)=\|u-v\|_{L^{\infty}\left((0,T), L^{2}(\Gamma)\right)}.
$$
By \cite[Theorem 1.2.5]{Caz03}, $(E_1, d_1)$ is a complete metric space.
 We now consider $\mathcal{H}$ defined by
\begin{equation*}\label{H}
\mathcal{H}(u)(t)=e^{i\Delta_\beta t} u_0+\mathcal{G}(u)(t),\qquad \mathcal{G}(u)(t)=i \int_{0}^{t} e^{i\Delta_\beta(t-s)} g(u(s)) d s
\end{equation*}
for all $u \in E_1$ and  $t \in (0, T) .$ We note that if $u \in L^{\infty}\left((0, T), H^1(\Gamma)\right),$ then $g(u) \in L^{\infty}\left((0, T), H^1(\Gamma)\right)$ by \eqref{g_1}.  Using \eqref{group} and \eqref{g_1}, for every $u \in E_1$ we obtain
$$
\|\mathcal{H}(u)(t)\|_{L^{\infty}\left((0, T), H^1(\Gamma)\right)} \leq C\|u_0\|_{H^1}+TC\|g(u)\|_{L^{\infty}\left((0, T), H^1(\Gamma)\right)} \leq C\|u_0\|_{H^1}+CTC_1(M) M.
$$
Furthermore, it follows from \eqref{g_1} that for  $u, v \in E_1$
$$
\|\mathcal{H}(u)(t)-\mathcal{H}(v)(t)\|_{L^{\infty}\left((0,T), L^{2}(\Gamma)\right)} \leq T C_2(M)\|u-v\|_{L^{\infty}\left((0,T), L^{2}(\Gamma)\right)}.
$$
Therefore, we see that if
$$
M=2C\|u_0\|_{H^1},  \qquad C T C_1(M) \leq \frac{1}{2}, \qquad T C_2(M) <1,
$$
then $\mathcal{H}$ is a strict contraction of $(E_1, d_1)$. Thus, it has a fixed point $u(t)$ which is a solution to \eqref{NLS}. 
Finally,  let $t, t_n\in [0, T]$ and $t_n\underset{n\to\infty}{\longrightarrow}t$, then

\begin{align*}\|\mc G(u)(t)-\mc G(u)(t_n)\|_{H^1}\leq \int_0^{t}\|\Big(e^{i\Delta_\beta(t-s)}-e^{i\Delta_\beta(t_n-s)}\Big)g(u(s))\|_{H^1}ds+\left|\int_t^{t_n}\|e^{i\Delta_\beta(t_n-s)}g(u(s))\|_{H^1}ds\right|.\end{align*}
  By continuity of $e^{i\Delta_\beta t}$ in $H^1(\Gamma)$, we conclude   $\mathcal{G}(u) \in C\left([0,T], H^1(\Gamma)\right),$  therefore $u(t)\in C\left([0,T], H^1(\Gamma)\right).$

 Uniqueness follows by the Gronwall lemma. Indeed, let $u_1, u_2$ be two solutions of \eqref{NLS} and \\ $\widetilde M=\sup\limits_{t\in[0,T]}\max\{\|u_1(t)\|_{H^1}, \|u_2(t)\|_{H^1}\}$, then, by \eqref{g_1}, one gets 
  $$\|u_1(t)-u_2(t)\|_2\leq C_2(\widetilde M)\int_0^t\|u_1(s)-u_2(s)\|_2ds.$$  
   The proof of the  blow-up alternative follows standardly  by a bootstrap argument.

\noindent\textit{ Step 3.} We show continuous dependence. Let $u_0 \in H^1(\Gamma)$ and consider $\{u_0^n\}_{n\in \mathbb{N}} \subset H^1(\Gamma) $ such that $u^{n}_0 \underset{n\to \infty}{\longrightarrow} u_0$ in $H^1(\Gamma)$. Assume that $u^{n}$ is  the maximal solution corresponding to the initial value $u_0^n$.

Since $\left\|u_0^n\right\|_{H^1} \leq 2\|u_0\|_{H^1}$ for $n$ sufficiently large, we deduce from \textit{Step 2} that there exists $T=T\left(\|u_0\|_{H^1}\right)$ such that $u$ and $u^{n}$ are defined on $[0, T]$ for $n \geq n_{0}$, and
$$
\|u\|_{L^{\infty}\left((0, T), H^1(\Gamma)\right)}+\sup _{n \geq n_{0}}\left\|u^{n}\right\|_{L^{\infty}\left((0, T), H^1(\Gamma)\right)} \leq 6C\|u_0\|_{H^1}.
$$
The second estimate in \eqref{g_1} yields  that there exists $C_3=C_3(\|u_0\|_{H^1})$ such that
$$
\left\|u(t)-u^{n}(t)\right\|_{2} \leq\left\|u_0-u^{n}_0\right\|_{2}+C_3 \int_{0}^{t}\left\|u(s)-u^{n}(s)\right\|_{2} d s.
$$
 We then conclude,  by the Gronwall lemma,
\begin{equation}\label{cont_dep} 
\left\|u(t)-u^{n}(t)\right\|_2 \leq\left\|u_0-u^{n}_0\right\|_{2} e^{T C_3} \underset{n \rightarrow \infty}{\longrightarrow} 0.
\end{equation}
Therefore,  from  \eqref{g_2} and \eqref{cont_dep} it follows that there exist $C_4=C_4(\|u_0\|_{H^1})$ and $\varepsilon_{n} \downarrow 0$ such that
$$
\left\|u(t)-u^{n}(t)\right\|_{H^1} \leq \varepsilon_{n}+C\left\|u_0-u^{n}_0\right\|_{H^1}+C_4  \int_{0}^{t}\left\|u(s)-u^{n}(s)\right\|_{H^1} d s
$$
for all $t \in [0, T]$. Applying again Gronwall's lemma, we get the result.

\noindent\textit{Step 4. } We show that for $p>2$ the mapping data-solution is of class $C^2$. 
 
 We interpret the  solution   $u(t)$ to \eqref{NLS} as the fixed point of the mapping $\mathcal H(u)(t)$ defined above in the metric space $(E_2,d_2)$: 
 $$
E_2=\left\{u \in C\left([0,T], H^1(\Gamma)\right):\|u\|_{C\left([0,T], H^{1}(\Gamma)\right)} \leq M\right\}, \quad {d}_2(u, v)=\|u-v\|_{C\left([0,T], H^1(\Gamma)\right)}.
$$
Observe that $g(u)\in C^2(\mathbb{C}, \mathbb{C})$, which implies 
\begin{equation}\label{g_3}
\|g(u)-g(v)\|_{H^1}\leq C_3(M)\|u-v\|_{H^1} \quad\text{for}\quad \|u\|_{H^1}\leq M,\,\,\|v\|_{H^1}\leq M.
\end{equation}
By \eqref{g_1},  for every $u \in E_2$ 
$$
\|\mathcal{H}(u)(t)\|_{C\left([0, T], H^1(\Gamma)\right)} \leq C\|u_0\|_{H^1}+TC\|g(u)\|_{C\left([0, T], H^1(\Gamma)\right)} \leq C\|u_0\|_{H^1}+CTC_1(M) M.
$$
Furthermore, it follows from \eqref{g_3} that, if $u, v \in E_2,$ then
$$
\|\mathcal{H}(u)(t)-\mathcal{H}(v)(t)\|_{C\left([0, T], H^1(\Gamma)\right)} \leq C T C_3(M)\|u-v\|_{C\left([0, T], H^1(\Gamma)\right)}.
$$
Let
\begin{equation}\label{estim_contract}
M=2C\|u_0\|_{H^1},  \qquad C T C_1(M) \leq \frac{1}{2}, \qquad C T C_3(M) < 1,
\end{equation}
then $\mathcal{H}$ is a contraction of $(E_2, d_2)$.

Consider the mapping
\begin{align*}
&\mathcal{J}: B(u_0, u(t))\subset H^1(\Gamma)\times C\left([0, T], H^1(\Gamma)\right)\longrightarrow C\left([0, T], H^1(\Gamma)\right),\\
&\mathcal{J}(v_0, v(t))=v(t)-{e}^{i \Delta_{\beta}t} v_{0}-i \int_{0}^{t} e^{i \Delta_{\beta}(t-s)} g(v(s)) d s.
\end{align*}
Here $B(u_0, u(t))$ is an open neighborhood of  $(u_0, u(t)).$ It is obvious that $\mathcal{J}(u_0, u(t))=0$.
One has 
$$D_{v(t)}\mathcal{J}(u_0, u(t))h(t)=h(t)- i \int_{0}^{t} e^{i \Delta_{\beta}(t-s)}\left[\partial_{v}g(u)h+\partial_{\overline{v}}g(u)\overline{h}\right](s)ds=(I-A)h(t),$$ with
$Ah(t)=i \int_{0}^{t} e^{i \Delta_{\beta}(t-s)}\left[\partial_{v}g(u)h+\partial_{\overline{v}}g(u)\overline{h}\right](s)ds$. From $CTC_1(M)\leq \frac{1}{2}$ we conclude 
$$\|Ah(t)\|_{C\left([0, T], H^1(\Gamma)\right)}\leq  \frac{1}{2}\|h(t)\|_{C\left([0, T], H^1(\Gamma)\right)},$$ hence $D_{v(t)}\mathcal{J}(u_0, u(t))$  is invertible  on  $C\left([0, T], H^1(\Gamma)\right)$ (as a linear operator on a real Banach space), and consequently it is a bijection. Therefore, by the Implicit Function Theorem, we conclude the existence of an open neighborhood $B(u_0)$ of $u_0$ and a unique function $f: B(u_0)\longrightarrow C\left([0, T], H^1(\Gamma)\right)$ such that   $\mathcal{J}(v_0, f(v_0))=0$ for all $v_0\in B(u_0)$, i.e. 
$$f(v_0)={e}^{i \Delta_{\beta}t} v_{0}+i \int_{0}^{t} e^{i \Delta_{\beta}(t-s)} g(f(v_0))(s) d s.$$ 
Hence $f(v_0)$ is the solution to \eqref{NLS} corresponding to the initial value  $v_0$.
Finally, since $g(u)$ is of class $C^2$ for $p>2$, we get that $f$ is $C^2$ mapping.  

The  proof of   conservation laws \eqref{conserv}  might be obtained involving the  regularization procedure analogous to the one introduced in the proof of \cite[Theorem 3.3.5]{Caz03}. 
\end{proof}
\begin{remark}\label{rem_glob_H_1}
For $1<p<5$ problem \eqref{NLS} is globally well-posed in $H^{1}(\Gamma)$.  This follows from \cite[Theorem 3.4.1]{Caz03} (in particular, see formula $(3.4.1)$). Namely, it is sufficient to observe  
\begin{equation}\label{glob_rem}
\begin{array}{l}
\|u\|_{p+1}^{p+1}-\frac{1}{\beta}\left|\sum\limits_{j=1}^N u_{j}(0)\right|^{2} \leq C_1\|u^{\prime}\|_{2}^{\frac{p-1}{2}}\|u\|_{2}^{\frac{p+3}{2}}+ \varepsilon\|u^{\prime}\|_{2}^{2}+C_{2}\|u\|_{2}^{2} \\
\leq  2\varepsilon\|u^{\prime}\|_{2}^{2}+C_{3}\|u\|_{2}^{\frac{2(p+3)}{5-p}}+C_{2}\|u\|_{2}^{2} \leq 2 \varepsilon\|u\|_{H^{1}}^{2}+C\left(\left\|u_{0}\right\|_{2}\right).
\end{array}
\end{equation}
The above estimate follows from the conservation of charge, estimate
\begin{equation}\label{v(0)}
\frac{1}{|\beta|}\left|\sum\limits_{j=1}^N v_{j}(0)\right|^{2}\leq \left|\frac{N^2}{\beta}\right|\|v\|_{\infty}^{2} \leq C\|v^{\prime}\|_{2}\|v\|_{2} \leq \varepsilon\|v^{\prime}\|_{2}^{2}+C_{\varepsilon}\|v\|_{2}^{2}, 
\end{equation}
  the Gagliardo-Nirenberg inequality for $v\in H^1(\Gamma)$
$$
\|v\|_{r} \leq C\|v'\|^{\alpha}_2\|v\|_{q}^{1-\alpha}, \quad r, q \in[2,+\infty], \quad  r \geq q, \quad \alpha=\frac{2}{2+q}(1-q / r), 
$$
  and the Young inequality  $ab\leq \delta a^{q}+C_{\delta} b^{q^{\prime}}, \frac{1}{q}+\frac{1}{q^{\prime}}=$ 1, $q, q^{\prime}>1, a, b \geq 0$. Observe that the key point is that $q=\frac{4}{p-1}>1$ for $1<p<5$.
\end{remark}
Below we prove the regularity of solution to \eqref{NLS} when $u_0\in \dom(\Delta_\beta).$
 Given the quantity
$$
0<m:=1-2 \inf \sigma(-\Delta_\beta)<\infty,
$$
we introduce the norm $\|v\|_{\beta}:=\|(-\Delta_\beta+m) v\|_{2}$ that endows $\operatorname{dom}(\Delta_\beta)$ with the structure of a Hilbert space. Observe that this norm for any real $\beta$ is equivalent to $H^{2}$-norm on the graph. Indeed, 
$$
\|v\|_{\beta}^{2}=\|v''\|_{2}^{2}+m^{2}\|v\|_{2}^{2}+2 m\|v'\|_{2}^{2}+ \frac{2 m}{\beta}\left|\sum\limits_{j=1}^N v_{j}(0)\right|^{2}.
$$
Due to the choice of $m$ and inequality \eqref{v(0)}, we get
$$
C_{1}\|v\|_{H^{2}(\Gamma)}^{2} \leq\|v''\|_{2}^{2}+m\|v\|_{2}^{2} \leq\|v\|_{\beta}^{2} \leq C_{2}\|v\|_{H^{2}(\Gamma)}^{2}.
$$
In what follows we will use the notation $D_{\beta}=\left(\operatorname{dom}(\Delta_\beta),\|\cdot\|_{\beta}\right)$.
\begin{proposition}\label{well_D_H}  Let $p>2$ and $u_0\in\dom(\Delta_\beta)$. Then there exists $T>0$  such that problem \eqref{NLS} has  a unique strong solution  $u(t)\in C\left([0,T], D_{\beta}\right)\cap C^1\left([0,T], L^2(\Gamma)\right)$.  Moreover, problem \eqref{NLS} has a maximal solution defined on
an interval of the form  $[0, T_{\beta})$, and the following ``blow-up alternative'' holds: either $T_{\beta} = \infty$ or $T_{\beta} <\infty$ and
$$\lim\limits_{t\to T_{\beta}}\|u(t)\|_{\beta} =\infty.$$
\end{proposition}
\begin{proof}
The proof is analogous to the one of  \cite[Theorem 2.3]{NataMasa2020} (for the NLS-$\delta$ equation with $p>4$) with a few modifications concerning the restriction $p>2$. 

 Let $T>0$ to be chosen later. We will use the notation $$X_\beta=C\left([0,T], D_\beta\right)\cap C^1\left([0,T], L^2(\Gamma)\right),$$  and equip  the space $X_\beta$ with the norm
$$\|u(t)\|_{X_\beta}=\sup\limits_{t\in [0,T]}\|u(t)\|_{\beta}+\sup\limits_{t\in [0,T]}\|\partial_t u(t)\|_2.$$
Consider $$E=\Big\{u(t)\in X_\beta:\, u(0)=u_0,\,\, \|u(t)\|_{X_\beta}\leq M \Big\},$$
where $M$ is a positive constant to be chosen later as well. 
It is obvious that  $(E, d)$ is a complete metric space  with the metric $d(u, v)=\|u-v\|_\beta$. 

Recall $$\mathcal{H}(u)(t)=e^{i\Delta_\beta t} u_0+\mathcal{G}(u)(t)=e^{i\Delta_\beta t} u_0+i \int_{0}^{t} e^{i\Delta_\beta(t-s)} g(u(s)) d s,\qquad g(u)=|u|^{p-1}u.
$$

\noindent\textit{Step 1.}  We will show that $\mathcal H:\, E\to X_\beta$.

 {\bf 1.}\,  It is known  that $\dom(\Delta_\beta)=\{v\in L^2(\Gamma):\, \lim\limits_{h\to 0}h^{-1}(e^{i\Delta_\beta h}-I)v\,\,\text{exists} \}.$ Obviously $w(t):=e^{i\Delta_\beta t}u_0\in \dom(\Delta_\beta)$ and $w(t)\in C\left([0, T], D_\beta\right)$.  Moreover, $\partial_t w(t)=i\Delta_\beta e^{i\Delta_\beta t}u_0=i\Delta_\beta w(t)$, then   $\partial_t w(t)\in C\left([0,T], L^2(\Gamma)\right).$ 
 
 {\bf 2.}\, The inclusion  $\mathcal G(u)(t)\in C^1\left([0,T],L^2(\Gamma)\right)$ follows rapidly.
Indeed, \cite[Lemma 4.8.4]{Caz03} implies that $\partial_tg(u(t))\in L^1\left((0,T), L^2(\Gamma)\right)$, and the formula  
 \begin{equation}\label{well_2a}
 \partial_t\mathcal G(u)(t)=ie^{i\Delta_\beta t}g(u(0))+i\int_0^te^{i\Delta_\beta(t-s)}\partial_sg(u(s))ds 
 \end{equation} from the proof of \cite[Lemma 4.8.5]{Caz03} induces $\mathcal G(u)(t)\in C^1\left([0,T],L^2(\Gamma)\right).$
 
{\bf 3.} Below we will show that $\mathcal G(u)(t)\in C\left([0,T], D_\beta\right)$.   First, we need to prove that $\mc G(u)(t)\in\dom(\Delta_\beta).$  
 Second inequality in \eqref{g_1} implies $g(u(t))\in C\left([0,T],L^2(\Gamma)\right)$, then for $t\in [0, T)$ and $h\in(0, T-t]$ we get
 \begin{equation}\label{well_1}
 \begin{split} &\frac{e^{i\Delta_\beta h} -I}{h}\mc G(u)(t)=\frac{1}{h}\int_0^t e^{i\Delta_\beta (t+h-s)}g(u(s))ds-\frac{1}{h}\int_0^te^{i\Delta_\beta (t-s)}g(u(s))ds\\&=\frac{\mc G(u)(t+h)-\mc G(u)(t)}{h}-\frac{1}{h}\int_t^{t+h}e^{i\Delta_\beta (t+h-s)}g(u(s))ds.
 \end{split}
 \end{equation}
 Letting $h\to 0$, by the Mean Value Theorem, we arrive at  $i\Delta_\beta\mc G(u)(t)=\mc G(u)'(t)-g(t)$, i.e.  we obtain the existence of the limit in \eqref{well_1},   therefore, $\mc G(u)(t)\in \dom(\Delta_\beta)$.    This is still true for $t=T$ since operator $-\Delta_\beta$ is closed.
 Note that we have used differentiability of $\mc G(u)(t)$ proved above.

It remains to prove the continuity of $\mc G(u)(t)$ in $\beta$-norm.
  We will use the integration by parts formula (see \cite[Proposition A.1]{NataMasa2020} with $H=-\Delta_\beta$)
 \begin{equation}\label{by_parts}
 \begin{split}
 &\mc G(u)(t)=\int_0^te^{i\Delta_\beta(t-s)}g(u(s))ds=-i(-\Delta_\beta+m)^{-1}g(u(t))+ie^{i\Delta_\beta t}(-\Delta_\beta+m)^{-1}g(u(0))\\&+m(-\Delta_\beta+m)^{-1}\int_0^te^{i\Delta_\beta(t-s)}g(u(s))ds+i(-\Delta_\beta+m)^{-1}\int_0^te^{i\Delta_\beta(t-s)}\partial_sg(u(s))ds.
 \end{split}
\end{equation} 
It is easily seen that 
\begin{equation}\label{partial_g}
\partial_sg(u(s))=\partial_u g(u)\partial_s u+\partial_{\overline u} g(u)\partial_s \overline{u}.
\end{equation}
Let $t_n, t\in [0,T]$  and $t_n\to t$. By \eqref{by_parts} and \eqref{partial_g}, we deduce
\begin{equation}\label{well_3}
\begin{split}
&\|\mc G(u)(t)-\mc G(u)(t_n)\|_\beta\leq \|g(u(t))-g(u(t_n))\|_2+m\int_0^{t}\|\Big(e^{i\Delta_\beta(t-s)}-e^{i\Delta_\beta(t_n-s)}\Big)g(u(s))\|_2ds\\&+m\left|\int_t^{t_n}\|e^{i\Delta_\beta(t_n-s)}g(u(s))\|_2ds\right|\\&+\int_0^{t}\|\Big(e^{i\Delta_\beta(t-s)}-e^{i\Delta_\beta(t_n-s)}\Big)\bigl(\partial_u g(u)\partial_su\bigr)(s)\|_2ds+\left|\int_t^{t_n}\|e^{i\Delta_\beta(t_n-s)}\bigl(\partial_u g(u)\partial_su\bigr)(s)\|_2ds\right|\\&+\int_0^{t}\|\Big(e^{i\Delta_\beta(t-s)}-e^{i\Delta_\beta(t_n-s)}\Big)\bigl(\partial_{\overline{u}} g(u)\partial_s\overline{u}\bigr)(s)\|_2ds+\left|\int_t^{t_n}\|e^{i\Delta_\beta(t_n-s)}\bigl(\partial_{\overline{u}} g(u)\partial_s\overline{u}\bigr)(s)\|_2ds\right|.
\end{split}
\end{equation}
Therefore, using \eqref{g_1},\eqref{well_3}, unitarity and continuity properties of $e^{i\Delta_\beta t}$,
we obtain continuity of $\mc G(u)(t)$ in $D_\beta$.

\noindent \textit{Step 2.} Now our aim is to choose  $T$ in order to guarantee invariance  of $E$ for the mapping $\mc H$, i.e.  $\mc H: E\to E$.

{\bf 1.} By \eqref{partial_g}, one has
\begin{equation}\label{well_4}
g(u(t))=\int_0^t\bigl(\partial_u g(u)\partial_su+\partial_{\overline{u}} g(u)\partial_s\overline{u}\bigr)(s)ds+g(u(0)).
\end{equation}
Let $u(t)\in E$ and $t\in [0,T]$. Using the Sobolev embedding, \eqref{by_parts}, \eqref{well_4},   and equivalence of $\beta$- and $H^2$-norm, we obtain
\begin{equation}\label{well_5}
\begin{split}
&\|\mc H(u)(t)\|_\beta\leq \|e^{i\Delta_\beta t}u_0\|_\beta+\|\int_0^te^{i\Delta_\beta(t-s)}g(u(s))ds\|_\beta\leq  \|u_0\|_\beta\\&+ \|\int_0^t\bigl(\partial_u g(u)\partial_su+\partial_{\overline{u}} g(u)\partial_s\overline{u}\bigr)(s)ds+g(u(0))\|_2\\&+\|g(u(0))\|_2+m\int_0^t\|g(u(s))\|_2ds+\int_0^t\|(\partial_ug(u)\partial_s u)(s)\|_2ds+\int_0^t\|(\partial_{\overline{u}}g(u)\partial_s \overline{u})(s)\|_2ds\\&\leq \|u_0\|_\beta+C_1\|u_0\|^p_\beta+C_2\int_0^t\|u\|_\infty^{p-1}\|\partial_su(s)\|_2ds+C_3\int_0^t\|u\|_\infty^{p-1}\|u(s)\|_2ds\\&\leq \|u_0\|_\beta+C_1\|u_0\|^p_\beta+K_1(M)TM^p.
\end{split}
\end{equation}
{\bf 2.} Below  we will estimate $\|\partial_t\mc H(u)(t)\|_2$.
Observe that 
\begin{equation}\label{well_6}
\|\partial_te^{i\Delta_\beta t}u_0\|_2=\|i\Delta_\beta u_0\|_2\leq \|u_0\|_\beta.
\end{equation}
Using  \eqref{well_2a}, \eqref{well_6},  we obtain the estimate
\begin{equation}\label{well_7}
\begin{split}
&\|\partial_t\mc H(u)(t)\|_2\leq \|u_0\|_\beta+\|g(u(0))\|_2+\int_0^t\|(\partial_ug(u)\partial_s u)(s)\|_2ds+\int_0^t\|(\partial_{\overline{u}}g(u)\partial_s \overline{u})(s)\|_2ds\\&\leq \|u_0\|_\beta+C_1\|u_0\|^p_\beta+C_2\int_0^t\|u\|_\infty^{p-1}\|\partial_s u(s)\|_2ds\leq \|u_0\|_\beta+C_1\|u_0\|^p_\beta+K_2(M)TM^p.
\end{split}
\end{equation}
Finally, combining \eqref{well_5} and \eqref{well_7}, we arrive at
$$\|\mc H(u)(t)\|_{X_\beta}\leq 2\|u_0\|_\beta+2C_1\|u_0\|^p_\beta+(K_1(M)+K_2(M))TM^p.$$
We now let
$$\frac{M}{2}=\left(2\|u_0\|_\beta+2C_1\|u_0\|^p_\beta\right).$$
By choosing  $T\leq\dfrac{1}{2(K_1(M)+K_2(M))M^{p-1}}$, we get 
$$\|\mc H(u)(t)\|_{X_\beta}\leq M,$$  therefore, $\mc H: E\to E$.

\noindent\textit{Step 3.} Now we will choose $T$ to guarantee that $\mc H$ is a strict contraction on $(E,d)$. 
First,  since $g(u)$ is of class $C^2$, we get 
\begin{equation}\label{well_8}
|\partial_ug(v)-\partial_ug(w)|\leq K_3(M)|v-w|, \quad  |\partial_{\overline u}g(v)-\partial_{\overline u}g(w)|\leq K_3(M)|v-w|
\end{equation}
for $|v|, |w|\leq M$.

Let $v, w\in E$. From \eqref{g_1},\eqref{by_parts}, \eqref{well_4}, \eqref{well_8} it follows that 
\begin{equation}\label{well_9}
\begin{split}
&\|\mc H(v)(t)-\mc H(w)(t)\|_\beta=\|\int_0^te^{i\Delta_\beta(t-s)}\Big[g(v(s))-g(w(s))\Big]ds\|_\beta\\&\leq
m\int_0^t\|g(v(s))-g(w(s))\|_2ds+2\int_0^t\|\big(\partial_ug(v)\partial_sv-\partial_ug(w)\partial_s w\big)(s)\|_2ds\\&+2\int_0^t\|\big(\partial_{\overline u}g(v)\partial_s\overline v-\partial_{\overline u}g(w)\partial_s \overline w\big)(s)\|_2ds\leq m C_2(M)T\|v-w\|_2\\&+2\int_0^t\|\big(\partial_ug(v)-\partial_ug(w)\big)\partial_sv(s)\|_2ds+2\int_0^t\|\partial_u g(w)\big(\partial_sv-\partial_s w\big)(s)\|_2ds\\&+ 2\int_0^t\|\big(\partial_{\overline u}g(v)-\partial_{\overline u}g(w)\big)\partial_s\overline v(s)\|_2ds+2\int_0^t\|\partial_{\overline u} g(w)  \big(\partial_s\overline v-\partial_s \overline w\big)(s)\|_2ds \\& \leq m C_2(M)T\|v-w\|_2+4K_3(M)T\|\partial_sv\|_{\infty}\|v-w\|_2+4Tp\|w\|^{p-1}_\infty\|\partial_sv-\partial_s w\|_2\\&\leq
TK_4(M)\|v-w\|_{X_\beta}.
\end{split}
\end{equation}
{\bf 2.} To get the contraction property of $\mc H$, we need to estimate $L^2$-part of $X_\beta$-norm of $\mc H(v)(t)-\mc H(w)(t)$.
From \eqref{well_2a} we deduce
\begin{equation}\label{well_15}
\|\partial_t\mc H(v)(t)-\partial_t\mc H(w)(t)\|_2\leq \int_0^t\|\partial_sg(v(s))-\partial_sg(w(s))\|_2ds.
\end{equation}
Using \eqref{partial_g},\eqref{well_8}, from \eqref{well_15}
we get 
\begin{equation}\label{well_16}
\|\partial_t\mc H(v)(t)-\partial_t\mc H(w)(t)\|_2\leq K_5(M)T\|v-w\|_{X_\beta},
\end{equation} and finally from \eqref{well_9},\eqref{well_16} we  obtain
$$\|\mc H(v)(t)-\mc H(w)(t)\|_{X_\beta}\leq (K_4(M)+K_5(M))T\|v-w\|_{X_\beta}.
$$
Thus, for $$T<\min\left\{\frac{1}{2(K_1(M)+K _2(M))M^{p-1}}, \frac{1}{K_4(M)+K_5(M)}\right\}$$ 
the mapping $\mc H$ is the strict contraction of $(E,d)$. Therefore, by the  Banach fixed point theorem,  $\mc H$ has a unique fixed point $u\in E$ which is a solution of \eqref{NLS}.

Uniqueness of the solution follows in a standard way by  Gronwall's lemma.
The blow-up alternative can be shown by a bootstrap argument.
\end{proof}
\begin{remark}
For $1<p<5$, problem \eqref{NLS} is globally well-posed in $D_\beta$. 
By formula $(3.4.2)$ in the proof of \cite[Theorem 3.4.1]{Caz03} (see also formula \eqref{glob_rem} in Remark \ref{rem_glob_H_1}),   we obtain that $\|u(t)\|_{H^1}$ is uniformly bounded. In particular, $\|u(t)\|_\infty\leq C(\|u_0\|_2, E(u_0)), t\in \mathbb{R}.$ Then  from \eqref{well_5}, \eqref{well_7} we get 
\begin{align*}&\|u(t)\|_{X_\beta}\leq 2\left\|u_{0}\right\|_{\beta}+2 C_{1}\left\|u_{0}\right\|_{\beta}^{p}+2C_{2} \int_{0}^{t}\|u\|_{\infty}^{p-1}\|\partial_{s} u(s)\|_{2} d s+C_{3} \int_{0}^{t}\|u\|_{\infty}^{p-1}\|u(s)\|_{2} d s\\
&\leq C_1(\left\|u_{0}\right\|_{\beta})+C_2(\|u_0\|_2, E(u_0)) \int_{0}^{t}\left\| u(s)\right\|_{X_\beta} d s, \quad t\in[0, T].
\end{align*}
By Gronwall's lemma,  we get
$$\|u(t)\|_{X_\beta}\leq  C_1(\left\|u_{0}\right\|_{\beta})e^{TC_2(\|u_0\|_2, E(u_0))}. $$
Therefore, using a standard bootstrap argument, one can show the global existence.  
\end{remark}
We end this section  mentioning that  analogously to \cite[Proposition 2.5]{NataMasa2020},  the following  ``virial identity" result holds.
\begin{lemma}\label{virial}
Assume that  $u_0\in \Sigma(\Gamma)=\{v\in H^1(\Gamma):\, xv\in L^2(\Gamma)\}$ and 
$u(t)$ is the corresponding maximal solution to \eqref{NLS}. 
Then $u(t)\in C\left([0, T_{H^1}), \Sigma(\Gamma)\right)$, and
the function 
\begin{equation*}\label{well_17}
f(t):=\int\limits_\Gamma x^2|u(t,x)|^2dx=\|xu(t)\|_2^2
\end{equation*}
belongs to $C^2[0,T_{H^1})$. Moreover,
\begin{equation}\label{well_17a}
f'(t)=4\im \int\limits_\Gamma x\overline{u}\partial_x u\,dx\quad\text{and }\quad
f''(t)=8 P(u(t)),\quad t\in [0, T_{H^1}),
\end{equation}
where 
\begin{align}\label{P}
P(v)=\|v'\|^{2}_2-\frac{1}{2\beta}\left|\sum\limits_{j=1}^Nv_{j}(0)\right|^{2}-\frac{p-1}{2(p+1)}\left\|v\right\|^{p+1}_{p+1}, \quad v\in H^1(\Gamma).
\end{align}
\end{lemma}
\section{Existence of the ground state}
By a standing wave of \eqref{NLS}, we mean a solution of the form $e^{i\omega t}\phi(x)$, where $\omega\in\mathbb{R}$ and $\phi$ is a solution of the stationary equation
\begin{equation}\label{S}
 -\Delta_\beta\phi+\omega\phi-\left|\phi\right|^{p-1}\phi=0.
\end{equation}
The above equation is the Euler-Lagrange equation of the action functional on $H^{1}(\Gamma)$
\begin{align*}
S_{\omega}(v):&=\frac{1}{2} {F}_{\beta}(v)+\frac{\omega}{2}\left\|v\right\|^{2}_2-\frac{1}{p+1}\|v\|^{p+1}_{p+1}.
\end{align*}
 The action $S_{\omega}$ is unbounded from below. Nevertheless, it is bounded from below on the so-called natural (or Nehari) manifold $\left\{v \in H^1(\Gamma): I_{\omega}(v)=0\right\}$, where $$I_{\omega}(v):= {F}_{\beta}(v)+\omega\left\|v\right\|^{2}_2-\|v\|^{p+1}_{p+1}.$$ Note that $I_{\omega}\left(v\right)=\langle S_{\omega}^{\prime}(v), v\rangle_{(H^1)'\times H^1}$, therefore, the Nehari manifold contains all the solutions to the stationary equation \eqref{S}. 
We consider the minimization problem on the Nehari manifold
\begin{align}\label{d}
d_{\omega}&=\inf\left\{S_{\omega}(v)\,:\,\,v\in H^{1}(\Gamma)\backslash \{0\},\,\,I_{\omega}(v)=0\right\}.
\end{align}
Our first result states the existence of the minimizer for $d_{\omega}$.
\begin{theorem}\label{GS}
Let $p>1$ and  $\omega>\frac{N^2}{\beta^2}. $ Then there exists $\beta^{\ast}<0$ such that problem \eqref{d}  admits a solution for any $\beta\in(\beta^*, -\frac{N}{\sqrt{\omega}})$.
\end{theorem}
\begin{proof}

\textit{Step 1.} Let 
\begin{align*}
&S_{\omega}^\infty(v)=\frac{1}{2}\|v'\|_2^2+\frac{\omega}{2}\left\|v\right\|^{2}_2-\frac{1}{p+1}\|v\|^{p+1}_{p+1},\quad
I_{\omega}^\infty(v)=\|v'\|_2^2+\omega\left\|v\right\|^{2}_2-\|v\|^{p+1}_{p+1},\\
&d^\infty_{\omega}=\inf\left\{S^\infty_{\omega}(v):\,\,v\in H^{1}(\Gamma)\backslash \{0\},\,\,I^\infty_{\omega}(v)=0\right\}= \inf\left\{\tfrac{p-1}{2(p+1)}\|v\|^{p+1}_{p+1}:\,\,v\in H^{1}(\Gamma)\backslash \{0\},\,\,I^\infty_{\omega}(v)=0\right\}.
\end{align*}
Index $\infty$ means that formally we assume  $\beta=\infty$.

Firstly, we show that  
\begin{equation}\label{d_inf_value}d^\infty_{\omega}=\tfrac{1}{2}\left(\tfrac{p+1}{2}\right)^{\frac{2}{p-1}} \omega^{\frac{p+3}{2(p-1)}} \int_{0}^{1}\left(1-t^{2}\right)^{\frac{2}{p-1}} d t.
\end{equation}
 Let $v\in H^1(\Gamma)$ and $v^*\in H^1_{\eq}(\Gamma)$ be its symmetric rearrangement (see \cite[Appendix A]{AQFF} for the definition). 
 We are going to substitute the minimizing problem on $H^1(\Gamma)$ by the one on $H^1_{\eq}(\Gamma)$. Notice that 
 \begin{equation*}\label{rearrange}
\|v\|_2=\left\|v^{*}\right\|_2,  \quad\|v\|_{p+1}=\left\|v^{*}\right\|_{p+1} \quad\|v^{\prime}\|_2\geq  \frac{2}{N}\|v^{* \prime}\|_2.
\end{equation*}
Since $I_\omega^\infty(v)\geq \frac{4}{N^{2}}\|v^{*\prime}\|^{2}_2-\|v^*\|_{p+1}^{p+1}+\omega\|v^*\|^{2}_2$, we get 
$$d_\omega^\infty\geq \inf \left\{\tfrac{p-1}{2(p+1)}\|v\|_{p+1}^{p+1}: v\in H^1_{\eq}(\Gamma)\backslash \{0\},\quad \tfrac{4}{N^{2}}\|v^{\prime}\|^{2}_2-\|v\|_{p+1}^{p+1}+\omega\|v\|^{2}_2 \leq 0 \right\}=:I_1.$$
 Taking the  rescaling $\lambda^{1/2} v(\lambda x)$ of $v(x)\in H^1_{\eq}(\Gamma)$ with  $\lambda=\left(\frac{N}{2}\right)^{\frac{4}{5-p}}$,
one arrives at 
\begin{align*}
I_1&=\left(\tfrac{N}{2}\right)^{\frac{2(p-1)}{5-p}}\inf \left\{\tfrac{p-1}{2(p+1)}\|v\|_{p+1}^{p+1}: v\in H^1_{\eq}(\Gamma)\backslash \{0\},\quad \left(\tfrac{N}{2}\right)^{\frac{2(p-1)}{5-p}}(\|v^{\prime}\|^{2}_2-\|v\|_{p+1}^{p+1})+\omega\|v\|^{2}_2 \leq 0 \right\}\\
&= \left(\tfrac{N}{2}\right)^{\frac{2(p-1)}{5-p}}\inf \left\{\tfrac{p-1}{2(p+1)}\|v\|_{p+1}^{p+1}: v\in H^1_{\eq}(\Gamma)\backslash \{0\},\quad \|v^{\prime}\|^{2}_2-\|v\|_{p+1}^{p+1}+\omega\left(\tfrac{2}{N}\right)^{\frac{2(p-1)}{5-p}}\|v\|^{2}_2 \leq 0 \right\}\\
 &
=N\left(\tfrac{N}{2}\right)^{\frac{2(p-1)}{5-p}} \inf \Big\{\tfrac{p-1}{2(p+1)}\|v\|_{L^{p+1}(\mathbb{R}^{+})}^{p+1}: v\in H^{1}(\mathbb{R}^{+})\setminus\{0\}, \\&\qquad\qquad\qquad\|v^{\prime}\|_{L^{2}(\mathbb{R}^{+})}^{2}-\|v\|_{L^{p+1}(\mathbb{R}^{+})}^{p+1}+\omega\left(\tfrac{2}{N}\right)^{{\frac{2(p-1)}{5-p}}}\|v\|_{L^{2}(\mathbb{R}^{+})}^{2} \leq 0\Big\}=:I_2.
\end{align*}
Let $\{v_n\}$ be a minimizing sequence for $I_2$. Then, taking even continuation of $v_n$ onto the whole line, we obtain 

\begin{equation}\label{ineq_long}
\begin{split}
 &I_2\geq \tfrac{N}{2}\left(\tfrac{N}{2}\right)^{\frac{2(p-1)}{5-p}} \inf \Big\{\tfrac{p-1}{2(p+1)}\|v\|_{L^{p+1}(\mathbb{R})}^{p+1}: v\in H^{1}(\mathbb{R}\setminus\{0\})\setminus\{0\}, \\&\qquad\qquad\qquad
\|v^{\prime}\|_{L^{2}(\mathbb{R})}^{2}-\|v\|_{L^{p+1}(\mathbb{R})}^{p+1}+\omega\left(\tfrac{2}{N}\right)^{{\frac{2(p-1)}{5-p}}}\|v\|_{L^{2}(\mathbb{R})}^{2} \leq 0\Big\}.
\end{split}
\end{equation}
By \cite[Lemma 4.4]{AdaNoj13}, 
\begin{equation}\label{d_line}
\begin{split}
&\inf \left\{\tfrac{p-1}{2(p+1)}\|v\|_{L^{p+1}(\mathbb{R})}^{p+1}: v\in H^{1}(\mathbb{R}\setminus\{0\})\setminus\{0\},\quad \|v^{\prime}\|_{L^{2}(\mathbb{R})}^{2}-\|v\|_{L^{p+1}(\mathbb{R})}^{p+1}+\omega\|v\|_{L^{2}(\mathbb{R})}^{2} \leq 0\right\}\\
&=\tfrac{p-1}{2(p+1)}\|\chi_{\mathbb{R}^+}\phi_\omega\|_{L^{p+1}(\mathbb{R})}^{p+1}=\frac{1}{2}\left(\tfrac{p+1}{2}\right)^{\frac{2}{p-1}} \omega^{\frac{p+3}{2(p-1)}} \int_{0}^{1}\left(1-t^{2}\right)^{\frac{2}{p-1}} d t,
\end{split}
\end{equation} 
where $\phi_\omega$
is the solution to the equation
$-\phi''(x)+\omega \phi(x)-|\phi(x)|^{p-1}\phi(x)=0,\, x\in\mathbb{R},$
given by
\begin{equation}\label{phi_omega}\phi_\omega(x)=\left\{\tfrac{(p+1)\omega}{2}\mathrm{sech}^2\left(\tfrac{(p-1)\sqrt{\omega}}{2}x\right)\right\}^{\frac{1}{p-1}}.
\end{equation}
Comparing the last line of \eqref{ineq_long} and \eqref{d_line}, we get (substituting $\omega$ by $\omega\left(\tfrac{2}{N}\right)^{{\frac{2(p-1)}{5-p}}}$)
\begin{equation*}
\begin{split}
&d^\infty_\omega\geq \tfrac{N}{2}\left(\tfrac{N}{2}\right)^{\frac{2(p-1)}{5-p}}\tfrac{1}{2}\left(\tfrac{p+1}{2}\right)^{\frac{2}{p-1}} \left[\omega\left(\tfrac{2}{N}\right)^{{\frac{2(p-1)}{5-p}}}\right]^{\frac{p+3}{2(p-1)}} \int_{0}^{1}\left(1-t^{2}\right)^{\frac{2}{p-1}} d t\\&=\tfrac{1}{2}\left(\tfrac{p+1}{2}\right)^{\frac{2}{p-1}} \omega^{\frac{p+3}{2(p-1)}} \int_{0}^{1}\left(1-t^{2}\right)^{\frac{2}{p-1}} d t.
\end{split}
\end{equation*}
Let $\tilde \phi(x)=(\tilde\phi_j(x))_{j=1}^N$ with
$
\tilde\phi_j(x)=\left\{\begin{array}{ll}
\phi_\omega,  & j=1 \\
0, & j \neq 1.
\end{array}\right.
$
It is easily seen that $I_\omega^\infty(\tilde \phi)=0$ and 
$$S_\omega^\infty(\tilde \phi)=\tfrac{1}{2}\left(\tfrac{p+1}{2}\right)^{\frac{2}{p-1}} \omega^{\frac{p+3}{2(p-1)}} \int_{0}^{1}\left(1-t^{2}\right)^{\frac{2}{p-1}} d t,$$
hence \eqref{d_inf_value} holds.

\noindent\textit{Step 2.} We prove that $d_{\omega}>0$.
 Since $\omega>\frac{N^2}{\beta^2}=-\inf\sigma(-\Delta_\beta)$,
 then  $F_{\beta}(v)+\omega\|v\|_{2}^{2}$ is equivalent to $H^1$-norm (see \cite[Lemma 4.13]{ArdCel21}). 
 Let $v \in H^{1}(\Gamma) \backslash\{0\}$ satisfy $I_{\omega}(v)=0 ,$ then
$$
\|v\|_{p+1}^{p+1}=F_{\beta}(v)+\omega\|v\|_{2}^{2}.
$$
 Summarizing the above, by the Sobolev embedding,  we obtain
$$\|v\|_{p+1}^{2} \leq C_{1}\|v\|_{H^{1}}^{2} \leq C_{2}\left(F_{\beta}(v)+\omega\|v\|_{2}^{2}\right)=C_{2}\|v\|_{p+1}^{p+1}, 
$$
therefore $C_{2}^{\frac{-1}{p-1}} \leq\|v\|_{p+1}$. Taking the infimum over $v$, we get $d_{\omega}>0$.

\noindent\textit{
 Step 3.}  We introduce  $\beta^*$ such that for $\beta^*<\beta$ one gets $d_\omega<d_\omega^\infty$. Consider $\phi_{\beta}=(\tilde\phi_{\beta})_{j=1}^N$, where 
\begin{equation}\label{phi_beta}\tilde\phi_{\beta}(x)=\left[\frac{(p+1) \omega}{2} \operatorname{sech}^{2}\left(\frac{(p-1) \sqrt{\omega}}{2} x-\tanh ^{-1}\left(\frac{N}{\beta \sqrt{\omega}}\right)\right)\right]^{\frac{1}{p-1}}.
\end{equation}
It is easy to check that $I_\omega(\phi_{\beta})=0$ (since $\phi_\beta$ satisfies \eqref{S}) and 
$$S_\omega(\phi_{\beta})=\tfrac{N}{2}\left(\tfrac{p+1}{2}\right)^{\frac{2}{p-1}} \omega^{\frac{p+3}{2(p-1)}} \int_{\frac{N}{|\beta|\sqrt{\omega}}}^{1}\left(1-t^{2}\right)^{\frac{2}{p-1}} d t,$$ then 
$$S_\omega(\phi_{\beta})<d_\omega^\infty\quad \Longleftrightarrow\quad  N\int_{\frac{N}{|\beta|\sqrt{\omega}}}^{1}\left(1-t^{2}\right)^{\frac{2}{p-1}} d t< \int_{0}^{1}\left(1-t^{2}\right)^{\frac{2}{p-1}} d t.
$$
Observe that $f(|\beta|)= N\int_{\frac{N}{|\beta|\sqrt{\omega}}}^{1}\left(1-t^{2}\right)^{\frac{2}{p-1}} d t$ is an increasing function on $(\frac{N}{\sqrt{\omega}}, \infty)$, and
$f((\frac{N}{\sqrt{\omega}}, \infty))=\left(0, N\int_{0}^{1}\left(1-t^{2}\right)^{\frac{2}{p-1}} d t \right)$. Therefore, there exists $\beta^*<-\frac{N}{\sqrt{\omega}}$ such that
\begin{equation*}\label{beta*}
N\int_{\frac{N}{|\beta^*|\sqrt{\omega}}}^{1}\left(1-t^{2}\right)^{\frac{2}{p-1}} d t= \int_{0}^{1}\left(1-t^{2}\right)^{\frac{2}{p-1}} d t.
\end{equation*}
Hence,
 for $\beta^*<\beta<-\frac{N}{\sqrt{\omega}}$ we get 
\begin{equation*}\label{comparison}
d_\omega\leq S_\omega(\phi_{\beta})<d_\omega^\infty.\end{equation*}
\textit{Step 4.}
Let
\begin{align*}
\tilde d_{\omega}:=\inf\left\{\tfrac{p-1}{2(p+1)}\left(F_\beta(v)+\omega\|v\|^2_2\right)\,:\,\,v\in H^{1}(\Gamma)\backslash \{0\},\,\,I_{\omega}(v)\leq0\right\}.
\end{align*}
We prove that $\tilde d_{\omega}=d_{\omega}$.
It is obvious that $\tilde d_\omega\leq d_\omega$.
Let $v\in H^{1}(\Gamma)\backslash \{0\}$ and $I_{\omega}(v)<0$. Put
\begin{align*}
\lambda_{1}:=\left(\frac{ {F}_{\beta}(v)+\omega\left\|v\right\|^2_2}{\left\|v\right\|^{p+1}_{{p+1}}}\right)^{\frac{1}{p-1}}.
\end{align*}
Then, since  $I_{\omega}(\lambda_{1}v)=0$ and $0<\lambda_{1}<1$ (this follows from the behavior of the function $g(\lambda)=I_\omega(\lambda v)$), we have
\begin{align*}
d_{\omega}\leq S_\omega(\lambda_1 v)=\tfrac{p-1}{2(p+1)}\left(F_\beta(\lambda_1 v)+\left\|\lambda_{1}v\right\|^{2}_{2}\right)=\tfrac{p-1}{2(p+1)}\lambda_{1}^{2}\left(F_\beta(v)+\left\|v\right\|^{2}_{2}\right)<\tfrac{p-1}{2(p+1)}\left(F_\beta(v)+\left\|v\right\|^{2}_{2}\right).
\end{align*}
Thus, we obtain $d_{\omega}\leq \tilde d_{\omega}$.

\noindent\textit{Step 5.}   Let   $\{v_{n}\}\subset H^{1}(\Gamma)\backslash \{0\}$ be a minimizing sequence for $d_{\omega}$. As long as
\begin{align}\label{mm1}
S_\omega(v_n)=\tfrac{p-1}{2(p+1)}\left( {F}_{\beta}(v_{n})+\omega\left\|v_{n}\right\|^2_2\right)=\tfrac{p-1}{2(p+1)}\left\|v_{n}\right\|^{p+1}_{p+1}\underset{n\to\infty}{\longrightarrow} d_{\omega},
\end{align}
 the sequence  $\{v_{n}\}$ is bounded in $H^{1}(\Gamma)$. Hence there exist a subsequence $\{v_{n_{k}}\}$ of $\{v_{n}\}$ and $v_{0}\in H^{1}(\Gamma)$ such that $\{v_{n_{k}}\}$ converges weakly to $v_{0}$ in $H^{1}(\Gamma)$. We may assume that  $v_{n_{k}}\neq 0$ and define
\begin{align*}
\lambda_{k}=\left(\frac{\|v'_{n_{k}}\|^2_2+\omega\left\|v_{n_{k}}\right\|^2_2}{\left\|v_{n_{k}}\right\|^{p+1}_{p+1}}\right)^{\frac{1}{p-1}}.
\end{align*}
Notice that $\lambda_{k}>0$ and $I^{\infty}_{\omega}(\lambda_{k}v_{n_{k}})=0$. Therefore, by\textit{ Step 3} and the definition of $d^{\infty}_{\omega}$, we obtain
\begin{align}\label{nn1}
d_{\omega}<d^{\infty}_{\omega}\leq\tfrac{p-1}{2(p+1)}\left\|\lambda_{k}v_{n_{k}}\right\|^{p+1}_{p+1}=\lambda^{p+1}_{k}\tfrac{p-1}{2(p+1)}\left\|v_{n_{k}}\right\|^{p+1}_{p+1}\,\,\,\mbox{for all}\,\,k\in \mathbb{N}.
\end{align} Furthermore, by  $I_\omega(v_{n_k})=0$, \eqref{mm1}, and the weak convergence, we get
\begin{align*}
\lim_{k\rightarrow\infty}\lambda_{k}&=\lim_{k\rightarrow\infty}\left(\frac{\left\|v_{n_{k}}\right\|^{p+1}_{p+1}-\frac{1}{\beta}\left|\sum\limits_{j=1}^N(v_{n_k})_j(0)\right|}{\left\|v_{n_{k}}\right\|^{p+1}_{p+1}}\right)^{\frac{1}{p-1}}=\left(\frac{d_{\omega}-\tfrac{p-1}{2(p+1)}\frac{1}{\beta}\left|\sum\limits_{j=1}^N(v_{0})_j(0)\right|}{d_{\omega}}\right)^{\frac{1}{p-1}}.
\end{align*}
 Taking the limit in \eqref{nn1}, we obtain $d_\omega<\lim\limits_{k\to \infty}\lambda_k^{p+1}d_\omega$.  Since, by \textit{Step 2}, $d_{\omega}>0$, we arrive at $\lim\limits_{k\to \infty}\lambda_k>1$, and consequently  $\frac{1}{\beta}\left|\sum\limits_{j=1}^N(v_{0})_j(0)\right|<0$. Thus, $v_{0}\neq0$.
  By the weak convergence, we obtain 
\begin{align}\label{Nop}
\lim_{k\rightarrow\infty}\left\{\left( {F}_{\beta}(v_{n_{k}})- {F}_{\beta}(v_{n_{k}}-v_{0})\right)+\omega\left(\left\|v_{n_{k}}\right\|^2_2-\left\|v_{n_{k}}-v_{0}\right\|^2_2\right)\right\} = {F}_{\beta}(v_{0})+\omega\left\|v_{0}\right\|^2_2.
\end{align}

%Next, passing to a subsequence of $\{v_{n_{k}}\}%$ if necessary, we may assume that $v_{n_{k}}%\underset{k\to\infty}{\longrightarrow} v_{0}$ %a.e. on $\Gamma$. 
  Therefore, by the Brezis-Leib lemma \cite{LBL},
$$\lim\limits_{k\rightarrow\infty}I_{\omega}(v_{n_{k}})-I_{\omega}(v_{n_{k}}-v_{0})=\lim\limits_{k\rightarrow\infty}-I_{\omega}(v_{n_{k}}-v_{0}) = I_{\omega}(v_{0}).$$

Since $v_{0}\neq0$, then the right-hand  side of \eqref{Nop} is positive. Hence it follows from \eqref{mm1} and \eqref{Nop} that
\begin{align*}
\tfrac{p-1}{2(p+1)}\lim_{k\rightarrow\infty}\left( {F}_{\beta}(v_{n_{k}}-v_{0})+\omega\left\|v_{n_{k}}-v_{0}\right\|^2_2\right)
<\tfrac{p-1}{2(p+1)}\lim_{k\rightarrow\infty}\left( {F}_{\beta}(v_{n_{k}})+\omega\left\|v_{n_{k}}\right\|^2_2\right)=d_{\omega}.
\end{align*}
Then, by \textit{Step 4} (using that $d_\omega=\tilde d_\omega$), we have $I_{\omega}(v_{n_k}-v_{0})>0$ for $k$ large enough. Thus, since $-I_{\omega}(v_{n_{k}}-v_{0})\underset{k\to\infty}{\longrightarrow} I_{\omega}(v_{0})$, we obtain $I_{\omega}(v_{0})\leq0$. By $d_\omega=\tilde d_\omega$ and the weak lower semicontinuity of norms, we conclude
\begin{align*}
d_{\omega}\leq\tfrac{p-1}{2(p+1)}\left( {F}_{\beta}(v_{0})+\omega\left\|v_{0}\right\|^2_2\right)\leq\tfrac{p-1}{2(p+1)}\lim_{k\rightarrow\infty}\left( {F}_{\beta}(v_{n_{k}})+\omega\left\|v_{n_{k}}\right\|^2_2\right)=d_{\omega}.
\end{align*}
Therefore, from \eqref{Nop}  we get
\begin{align*}
\lim_{k\rightarrow\infty} {F}_{\beta}(v_{n_{k}}-v_{0})+\omega\left\|v_{n_{k}}-v_{0}\right\|^2_2=0, 
\end{align*}
and, consequently,  we have $v_{n_{k}}\underset{k\to\infty}{\longrightarrow} v_{0}$ in $H^{1}(\Gamma)$ and $I_\omega(v_0)=0$.
\end{proof}
\begin{remark}
The restriction $\beta \in(\beta^{*},-\frac{N}{\sqrt{\omega}})$ seems to be non optimal at least for $\omega>\frac{p+1}{p-1}\frac{N^2}{\beta^2}$. For example, for even $N$    one gets $S(\phi_{\frac{N}{2}})<S(\phi_\beta)$, where $\phi_{\frac{N}{2}}$ is defined in Theorem \ref{phi_k} (the proof repeats the last  part of the proof of \cite[Theorem 5.3]{AdaNoj13}). Thus,     comparison of  $d_\omega^\infty$ with $S(\phi_{\frac{N}{2}})$ might guarantee the existence of the minimizer for some interval  larger than $(\beta^{*},-\frac{N}{\sqrt{\omega}}).$
\end{remark}
Above  we proved  the existence of the minimizer for sufficiently weak attractive  $\delta'_s  $ coupling. Nevertheless, below we show that the minimizer exists for any $\beta\in \mathbb{R}\setminus\{0\}$ and $\beta<0$ in two particular cases. The first case is related with $H^1_{\eq}(\Gamma),$ and the second case deals with  $H^1_{\frac{N}{2}}(\Gamma)$ (see two Lemmas below).
\begin{lemma}\label{lemma_eq}
Let $\beta\in \mathbb{R}\setminus\{0\}$. We set
\begin{align}\label{d^eq}
d_\omega^{\eq}=\inf\left\{S_\omega(v): v\in H^1_{\eq}(\Gamma)\setminus\{0\},\quad  I_\omega(v)=0\right\}.
\end{align}
Then  
$
d_\omega^{\eq}=S_\omega(\phi_{\beta}),$
where $\phi_{\beta}=(\tilde\phi_{\beta})_{j=1}^N$  is  defined by \eqref{phi_beta}.
\end{lemma}
\begin{proof}
Notice that 
$$d_\omega^{\eq}=Nd_\omega^{\mathrm{half}}=\frac{N}{2}d_{\omega, \mathrm{r}} ^{\mathrm{line}},$$
where  
     \begin{equation*}
       d_\omega^{\mathrm{half}}=\inf\left\{\begin{array}{c}\tfrac{1}{2}\|v'\|_{L^2(\mathbb{R}^+)}^2+\tfrac{\omega}{2}\|v\|_{L^2(\mathbb{R}^+)}^2-\frac{1}{p+1}\|v\|_{L^{p+1}(\mathbb{R}^+)}^{p+1}+\tfrac{N}{2\beta}|v(0)|^2:\\ \|v'\|_{L^2(\mathbb{R}^+)}^2+\omega\|v\|_{L^2(\mathbb{R}^+)}^2-\|v\|_{L^{p+1}(\mathbb{R}^+)}^{p+1}+\tfrac{N}{\beta}|v(0)|^2=0, \,v\in H^1(\mathbb{R}^+)\setminus\{0\}    
    \end{array}\right\}, 
        \end{equation*} 
 \begin{equation*}\label{well_27c}
    d_{\omega, \mathrm{r}}^{\mathrm{line}} =\inf\left\{\begin{array}{c}\tfrac{1}{2}\|v'\|_{L^2(\mathbb{R})}^2+\tfrac{\omega}{2}\|v\|_{L^2(\mathbb{R})}^2-\frac{1}{p+1}\|v\|_{L^{p+1}(\mathbb{R})}^{p+1}+\tfrac{N}{\beta}|v(0)|^2:\\ \|v'\|_{L^2(\mathbb{R})}^2+\omega\|v\|_{L^2(\mathbb{R})}^2-\|v\|_{L^{p+1}(\mathbb{R})}^{p+1}+\tfrac{2N}{\beta}|v(0)|^2=0, \,v\in H^1_{\mathrm{rad}}(\mathbb{R})\setminus\{0\}    
    \end{array}\right\}.
        \end{equation*}
         From  the results by \cite{RJJ, FukOht08} one gets  
    $$d_{\omega, \mathrm{r}} ^{\mathrm{line}}=   \tfrac{1}{2}\|\tilde\phi'_{\beta}\|_{L^2(\mathbb{R})}^2+\tfrac{\omega}{2}\|\tilde\phi_{\beta}\|_{L^2(\mathbb{R})}^2-\frac{1}{p+1}\|\tilde\phi_{\beta}\|_{L^{p+1}(\mathbb{R})}^{p+1}+\tfrac{N}{\beta}|\tilde\phi_{\beta}(0)|^2. $$
\end{proof}

\begin{lemma}\label{N/2}
Let $\beta<0$,  $N$ be even, and let the profile $\phi_{\frac{N}{2}}$ be defined in Theorem \ref{phi_k}. Set 
\begin{equation}\label{min_N_2}
d_{\omega}^{\frac{N}{2}}=\inf\Big\{S_{\omega}(v): v \in H_{\frac{N}{2}}^{1}(\Gamma) \setminus\{0\}, \quad I_{\omega}(v)=0\Big\}.
\end{equation}
Then 

$(i)$\, for  $\frac{N^2}{\beta^2}<\omega\leq \frac{p+1}{p-1}\frac{N^2}{\beta^2}$ the solution to \eqref{min_N_2} is given by $\phi_\beta$; 

$(ii)$\,  for  $\omega> \frac{p+1}{p-1}\frac{N^2}{\beta^2}$ the solution to \eqref{min_N_2} is given by $\phi_{\frac{N}{2}}$.
\end{lemma}
\noindent The proof of the above lemma follows from \cite[Theorem 4.1, Theorem 5.3]{AdaNoj13}, observing that for $v\in H^1_{\frac{N}{2}}(\Gamma)$ 
\begin{align*}
S_\omega(v)&=\frac{N}{4}\Big[\|v'_1\|^2_{L^2(\mathbb{R}^+)}+\|v'_N\|^2_{L^2(\mathbb{R}^+)}+\omega(\|v_1\|^2_{L^2(\mathbb{R}^+)}+\|v_N\|^2_{L^2(\mathbb{R}^+)})\\
&+\frac{N}{2\beta}\left|v_1(0)+v_N(0)\right|^2-\frac{2}{p+1}(\|v_1\|^{p+1}_{L^{p+1}(\mathbb{R}^+)}+\|v_N\|^{p+1}_{L^{p+1}(\mathbb{R}^+)})\Big].
\end{align*}
In other words, in $H^1_{\frac{N}{2}}(\Gamma)$ 
 minimizing problem \eqref{min_N_2} is ``equivalent'' to $\frac{N}{2}$ copies of the minimizing problem on the line.

\begin{remark}\label{V_remark}
In \cite{ArdCel21} the  NLS equation with the  $\delta$ coupling and decaying potential $V(x)$ on $\Gamma$
$$i \partial_{t} u(t)=-\Delta_{\gamma} u(t)+V(x) u(t)-|u(t)|^{p-1} u(t)$$
has been considered.  
Here  $\gamma<0$  and  $\left(-\Delta_{\gamma}v\right)(x)=(-v_j''(x))$ with  
\begin{equation*}\label{dom_delta}
\dom(\Delta_{\gamma})=\Big\{v\in H^1(\Gamma): v_1(0)=\ldots =v_N(0),\quad \sum^{N}_{j=1}v'_j(0)=\gamma v_{1}(0)\Big\}.
\end{equation*} 
It was shown that under the assumptions
\begin{equation}\label{V_cond}
\begin{split}
&V(x)=(V_j(x))_{j=1}^N \in L^{1}(\Gamma)+L^{\infty}(\Gamma),\qquad \lim _{x \rightarrow \infty} V_j(x)=0,\\& \int\limits_{\mathbb{R}^+} V_j(x)|\varphi(x)|^{2} d x<0,\,\,\varphi(x) \in H^{1}(\mathbb{R}^+)\setminus\{0\}
\end{split}\end{equation}
 there exist $0<\omega_0$ ($=-\inf\sigma(-\Delta_\gamma+V)$) and $\gamma^*<0$ such that for $\omega>\omega_0$ and $\gamma<\gamma^*$ problem \eqref{d}  admits a solution (in definition of   $S_\omega$ one needs to substitute  quadratic form $F_\beta$ by the quadratic form of the operator $-\Delta_\gamma+V$).

Under  assumptions  \eqref{V_cond}, one may consider   
\begin{equation}\label{V_beta}
i \partial_{t} u(t)=-\Delta_{\beta} u(t)+V(x) u(t)-|u(t)|^{p-1} u(t).
\end{equation}
It can be proven analogously that for $\omega>-\inf\sigma(-\Delta_\beta+V)$ and $\beta\in (\beta^*, -\frac{N}{\sqrt{\omega}})$ problem \eqref{d}  has a solution $\phi_{\beta, V}$ (with $S_\omega(v)$ substituted by $S_\omega(v)+\frac{1}{2}(Vv,v)_2$).
 \end{remark}
\section{Identification of the ground states}\label{4}
In this section we describe explicitly the solutions to \eqref{d}. In what follows we will use the notation 
\begin{equation}\label{phi_x}
\phi_{a}(x)=\left\{\tfrac{(p+1)\omega}{2}\mathrm{sech}^2\left(\tfrac{(p-1)\sqrt{\omega}}{2}(x+a)\right)\right\}^{\frac{1}{p-1}}, \quad a\in\mathbb{R}.
\end{equation} 
\begin{proposition}\label{ground}
Let $\omega>\frac{N^2}{\beta^2}$. Assume that  the   solution of  problem \eqref{d}  exists,   then it is given by the critical point of $S_\omega$. Moreover, each critical point of $S_\omega$ has the form
$\phi=e^{i\theta}(\phi_{x_j})_{j=1}^N,$
with $\theta\in \mathbb{R}$,
where $|x_j|=\frac{2}{(p-1)\sqrt{\omega}}\tanh^{-1}(t_j)$ and $t_j, \, j=1,\ldots,N,$ satisfy the system
\begin{equation}\label{system}
\left\{\begin{array}{c}
t_{1}^{p-1}-t_{1}^{p+1}=\ldots=t_N^{p-1}-t_{N}^{p+1}, \\
\sum\limits_{j=1}^N t_{j}^{-1}=|\beta|\sqrt{\omega}.
\end{array}\right.
\end{equation}
In particular, for $\beta>0$ ($\beta<0$) all the constants  $x_j$ are negative (positive). 
\end{proposition}
\begin{proof}
 Let $\phi$ be a minimizer. Since $I_{\omega}(\phi)=0$, we have
\begin{align}\label{lli}
\left\langle I^\prime_{\omega}(\phi),\phi\right\rangle_{(H^1)'\times H^1}=2\left( {F}_{\beta}(\phi)+\omega\left\|\phi\right\|^2_2\right)-(p+1)\left\|\phi\right\|^{p+1}_{p+1}=-(p-1)\left\|\phi\right\|^{p+1}_{p+1}<0.
\end{align} 
There exists a Lagrange multiplier $\mu\in\mathbb{R}$ such that $S^\prime_{\omega}(\phi)=\mu I^\prime_{\omega}(\phi)$. Furthermore, since
\begin{align*}
\mu\left\langle I^\prime_{\omega}(\phi),\phi\right\rangle_{(H^1)'\times H^1}=\left\langle S^\prime_{\omega}(\phi),\phi\right\rangle_{(H^1)'\times H^1}= I_{\omega}(\phi)=0,
\end{align*}
then, by \eqref{lli}, $\mu=0$. Hence $S^\prime_{\omega}(\phi)=0$.
As in the proof of \cite[Proposition 5.1]{AdaNoj13}, one can show that $S^\prime_{\omega}(\phi)=0$ is equivalent to \eqref{S} and $\phi\in \dom(\Delta_\beta)$.  The most general $L^2(\mathbb{R}^+)$-solution to 
$$-\varphi''+\omega \varphi-|\varphi|^{p-1} \varphi=0$$
is given by  $\sigma \phi_a(x), |\sigma|=1$ (see\eqref{phi_x}), therefore, $\phi=(e^{i\theta_j}\phi_{x_j})_{j=1}^N$.  Since $\phi$ is the minimizer of $S_\omega$,  one easily concludes $e^{i\theta_1}=\ldots=e^{i\theta_N}.$ From \eqref{dom_delta'} we get
\begin{equation}\label{system_cosh}\left\{\begin{array}{l}
\frac{\tanh \left(\frac{p-1}{2} \sqrt{\omega} x_{j}\right)}{\cosh ^{\frac{2}{p-1}}\left(\frac{p-1}{2} \sqrt{\omega} x_{j}\right)}-\frac{\tanh \left(\frac{p-1}{2} \sqrt{\omega} x_{j+1}\right)}{\cosh ^{\frac{2}{p-1}}\left(\frac{p-1}{2} \sqrt{\omega} x_{j+1}\right)}=0,\, j=1,\ldots, N-1, \\
\sum\limits_{j=1}^N\frac{1}{\cosh ^{\frac{2}{p-1}}\left(\frac{p-1}{2} \sqrt{\omega} x_{j}\right)}=-\beta\sqrt{\omega} \frac{\tanh \left(\frac{p-1}{2} \sqrt{\omega} x_{1}\right)}{\cosh ^{\frac{2}{p-1}}\left(\frac{p-1}{2} \sqrt{\omega} x_{1}\right)}.
\end{array}\right.\end{equation}
Introducing $t_j=\tanh \left(\frac{p-1}{2} \sqrt{\omega} |x_{j}|\right)$  and observing that $\cosh^{-2}\left(\frac{p-1}{2} \sqrt{\omega} x_{j}\right)=1-\tanh^2 \left(\frac{p-1}{2} \sqrt{\omega} x_{j}\right)$, we arrive at \eqref{system}. To complete the proof, notice that from the  first equation in \eqref{system_cosh} it follows that $x_j$ have the same sign.
\end{proof}
The next theorem gives a precise description of the family of the critical points of $S_\omega$ containing the ground state.
\begin{theorem}\label{phi_k}
Let $\omega>\frac{p+1}{p-1}\frac{N^2}{\beta^2}$, $N\geq 2$, $\beta\in\mathbb{R}\setminus\{0\}$. Then (up to permutation of the edges of the $\Gamma$ and rotation), the family of critical points of $S_\omega$ consists of $N$ profiles. The  first  $N-1$ profiles are given by $\phi_k=(\underset{1}{\phi_{x_1}},\ldots,\underset{k}{\phi_{x_1}},\underset{k+1}{\phi_{x_N}},\ldots,\underset{N}{\phi_{x_N}}),\, k=1,\ldots,N-1,$  
where $|x_1|=\frac{2}{(p-1)\sqrt{\omega}}\tanh^{-1}(t_1),$ $|x_N|=\frac{2}{(p-1)\sqrt{\omega}}\tanh^{-1}(t_N),$ and $0<t_1<t_N<1$ satisfy the system
\begin{equation}\label{system_k}
\left\{\begin{array}{c}
t_{1}^{p-1}-t_{1}^{p+1}=t_N^{p-1}-t_{N}^{p+1} \\
 kt_{1}^{-1}+(N-k)t_N^{-1}=|\beta|\sqrt{\omega}.
\end{array}\right.
\end{equation}
The $N$th critical point of $S_\omega$ is given by the symmetric profile  $\phi_\beta$. 
\end{theorem}
\begin{proof}
 Notice that  the first line of \eqref{system} might be rewritten  as 
\begin{equation*}\label{first_line}
f(t_1)=\ldots=f(t_N), \,\,\text{where}\,\,f(t):=t^{p-1}-t^{p+1}:[0,1)\rightarrow \mathbb{R}.
\end{equation*}
We have $f([0,1))=\left[0, \left(\frac{p-1}{2}\right)^{\frac{p-1}{2}}/\left(\frac{p+1}{2}\right)^{\frac{p+1}{2}}\right)$, and the maximum
$M=\left(\frac{p-1}{2}\right)^{\frac{p-1}{2}}/\left(\frac{p+1}{2}\right)^{\frac{p+1}{2}}$ is attained at the unique extremum point $t_{\max}=\sqrt{\frac{p-1}{p+1}}$. By monotonicity properties of $f$, we conclude that each $t_j$  may take only two possible values:   $t_1$,  $t_N$ (without loss of generality), where   $0<t_1<t_{\max}<t_N<1$ and $f(t_1)=f(t_N).$

Now, assuming that  the number of $t_j$ that take value $t_1$ (resp. $t_N$) is $k\in\{1,\ldots, N-1\}$ (resp. $N-k$), we get the second line of system \eqref{system_k}.

Then  expressing $t_1$ from the second line of \eqref{system_k}, by the first line, we get that $t_N$ satisfies  the equation $w(t_N)=0,$ where  
\begin{align*}
w(x)=\frac{x^2\left(a^{2}k^{p-1}-k^{p+1}\right)-2 a (N-k)k^{p-1}x+(N-k)^2k^{p-1}}{(ax-N+k)^{p+1}}+x^{2}-1, \quad a=|\beta|\sqrt{\omega}.
\end{align*}
Arguing as in the  proof  of \cite[Lemma 5.2]{AdaNoj13}, we can prove that  for any  $a>N\sqrt{\frac{p+1}{p-1}}$, there exists a unique root of $w(x)$ in  $\left(\frac{N}{a}, 1\right]$.  
It is sufficient to observe that:
\begin{itemize}

\item[$\bullet$]\, $w\left(\frac{N}{a}\right)=0$,

\item[$\bullet$]\, $
w'\left(\frac{N}{a}\right)=\frac{1}{ak}\left((p+1)N^2-a^2(p-1)\right)<0,$

\item[$\bullet$]\,$
w(1)=\frac{k^{p-1}}{(a-N+k)^{p+1}}((a-N+k)^2-k^2)>0,$
\item[$\bullet$]\,$
w'''(x)=\dfrac{-(p-1)p(p+1)k^{p-1}a^3(a^2-k^2)x^2+A_wx+C_w}{(ax-N+k)^{p+4}},$ where $A_w, C_w$ are real constants. Thus,  $w'''(x)<0$ for $x$  large enough.
\end{itemize}
To end the proof, we observe that assuming $k=N$, we get the solution $\phi_\beta$.

\end{proof}

\begin{remark}\label{k>N-k}
When   $N\geq 3$,   it seems that approach by \cite{AdaNoj13} cannot be applied to prove  that system \eqref{system_k} has a unique solution $\phi_\beta$ for $\frac{N^2}{\beta^2}<\omega\leq \frac{p+1}{p-1}\frac{N^2}{\beta^2}$. 
However, it still can be shown that the solution to \eqref{system_k} with $t_1<t_N$  for $k>N-k$  does not exist (it is sufficient to observe $g^2(t_1)<\left(\frac{k}{N-k}\right)^2g^2(t_N)$ in \cite[Formula (5.25)]{AdaNoj13}).  
 We conjecture that  for $\frac{N^2}{\beta^2}<\omega\leq \frac{p+1}{p-1}\frac{N^2}{\beta^2}$ the profile $\phi_\beta$ is the solution  to \eqref{d}.
\end{remark}
 \noindent We finish this section with the following important result which we will prove at the end of Subsection \ref{assym}.
\begin{theorem}\label{min_phi_1}
Let   $\beta<0, \omega>\frac{p+1}{p-1}\frac{N^2} {\beta^2}$, and assume that solution to problem \eqref{d} exists, then it is given by $\phi_1$.
\end{theorem} 
\section{Orbital and spectral instability}
The definition of the orbital stability involves the symmetry of the concrete Hamiltonian system in study. Since equation \eqref{NLS} is rotationally symmetric, we  
 define orbital stability as follows.
\begin{definition} The standing wave solution $u(t)=e^{i \omega t} \phi(x)$ is said to be orbitally stable in $H^1(\Gamma)$ by the flow of equation \eqref{NLS} if for any $\varepsilon>0$ there exists $\eta>0$ with the following property. If $u_{0} \in H^1(\Gamma)$ satisfies $\left\|u_{0}-\phi\right\|_{H^1}<\eta$, then the solution $u(t)$ of \eqref{NLS}  with $u(0)=u_{0}$ exists for any $t \in \mathbb{R}$ and
$$
\sup _{t \in \mathbb{R}} \inf _{\theta \in \mathbb{R}}\left\|u(t)-e^{i \theta} \phi\right\|_{H^1}<\varepsilon.
$$
Otherwise, the standing wave $u(t)=e^{i \omega t} \phi(x)$ is said to be orbitally unstable in $H^1(\Gamma)$.
\end{definition}
The definition of the\textit{ spectral instability} involves the concept of the linearization of \eqref{NLS} around the profile of the standing wave. After making  necessary technical steps we give  precise Definition \ref{linear_inst}.
\subsection{Instability analysis for symmetric profile $\phi_{\beta}$  via Grillakis/Jones approach}\label{symm_prof}
We begin this subsection with statement of one of the main results.
\begin{theorem}\label{main_instub}
Let $p>1$, $N\geq 2$, and $\phi_\beta$ be defined by \eqref{phi_beta}. If either $\beta<0$ and $\omega>\frac{p+1}{p-1}\frac{N^2}{\beta^2}$ or $\beta>0$ and $\frac{N^2}{\beta^2}<\omega<\frac{p+1}{p-1}\frac{N^2}{\beta^2}$,   then $e^{i\omega t}\phi_{\beta}$ is spectrally unstable. Moreover if, additionally, $p>2$, then   $e^{i\omega t}\phi_{\beta}$ is orbitally unstable. 
\end{theorem}
\begin{remark}
$(i)$ In the case $N=2$, the result of the above theorem was shown in \cite[Theorem 1.1]{AngGol20}  and \cite[Theorem 6.13]{AdaNoj13} (for $\beta>0$ and $\beta<0$, respectively). Moreover, our result extends \cite[Theorem 1.1]{AngGol20} for the case $p\in(3,5).$

$(ii)$ In \cite[Theorem 1.2]{AnGo2018} can be found some orbital stability/instability  results for $\beta<0$. In particular,  it was shown that for  $\frac{N^2}{\beta^2}<\omega<\frac{p+1}{p-1}\frac{N^2}{\beta^2}$ and $1<p\leq 5$ we have orbital stability of $e^{i\omega t}\phi_{\beta}$. Orbital instability   was shown for $\omega>\frac{p+1}{p-1}\frac{N^2}{\beta^2}, 1<p\leq 5 $, and  even $N$.

Moreover, for  $p>5$ and $\omega \neq \frac{p+1}{p-1}\frac{N^{2}}{\beta^{2}} $  we proved  that  there exists $\omega^{*}>\frac{N^{2}}{\beta^{2}}$ such that
$e^{i t \omega} \phi_{\beta}$ is orbitally unstable in $H^{1}(\Gamma)$ as $\omega>\omega^{*}$, and it is
orbitally stable in $H_{\eq}^{1}(\Gamma)$ as $\omega<\omega^{*}$.

\end{remark}

To define the concept of spectral instability, we linearize equation \eqref{NLS} around $\phi_{\beta}$. Firstly, we put $u(t)=e^{i\omega t}(\phi_{\beta}+v(t))$. Observing that $S_\omega$ is a $C^2$ functional,  $S'_\omega(\phi_{\beta})=0$, and   equation \eqref{NLS} has the form
$$i\partial_tu(t)=E'(u(t)),$$ we get
\begin{equation*}\label{linerize}
\partial_tv(t)=-iS''_\omega(\phi_{\beta})v(t)+O(\|v(t)\|^2_{H^1}).
\end{equation*}
  It is standard to verify that for $u\in H^1(\Gamma)$
\begin{equation*}\label{second_deriv}S''_\omega(\phi_{\beta})u=-\widetilde \Delta_\beta u+\omega u -(\phi_{\beta})^{p-1} u-(p-1)(\phi_{\beta})^{p-1}\re(u).\end{equation*}
Here the operator $-\widetilde \Delta_\beta$ is understood in the following sense: since the bilinear form \\ $t_\beta(u_1,  u_2)=(u'_1,  u'_2)_2+\frac{1}{\beta}\re\left(\sum\limits_{j=1}^Nu_{1,j}(0) \sum\limits_{j=1}^N\overline{u_{2,j}(0)}\right)$ is  bounded  on $H^1(\Gamma)$, there exists the unique bounded operator $\widetilde \Delta_\beta: H^1(\Gamma)\to (H^1(\Gamma))' $ such that $t_\beta(u_1, u_2)=\langle -\widetilde \Delta_\beta u_1, u_2\rangle_{(H^1)'\times H^1}.$
Notice that the bilinear  form 
\begin{equation*}\label{bil_form}
\begin{split}&b_\beta(u,v)=\langle S''_\omega(\phi_{\beta})u,v\rangle_{(H^1)'\times H^1}=\re\Big[\frac{1}{\beta}\sum\limits_{i,j=1}^Nu_{i}(0) \overline{v_{j}(0)}\\&+\int\limits_\Gamma\Big( u'\overline{v}'+\omega u\overline{v}-(\phi_{\beta})^{p-1} u\overline{v}-(p-1)(\phi_{\beta})^{p-1}\re(u)\overline{v}\Big)dx\Big]
\end{split}
\end{equation*}
is closed, densely defined, and bounded from below.
Then, by the Representation Theorem \cite[Chapter VI, Theorem 2.1]{Kat95}), 
we can associate with $b_\beta$ self-adjoint in $L^2(\Gamma)$ (with the real scalar product) operator
\begin{equation}\label{L_b}
\begin{split}
&L^\beta u=-\Delta_\beta u+\omega u-(\phi_{\beta})^{p-1} u-(p-1)(\phi_{\beta})^{p-1}\re(u),\quad \dom(L^\beta)=\dom(\Delta_\beta).
\end{split}
\end{equation}
We define a natural injection $\mathfrak{J}\in B\left(H^1(\Gamma),(H^1(\Gamma))'\right)$ by
$$\langle\mathfrak{J}u,v\rangle_{(H^1)'\times H^1}=(u,v)_2, \quad \text{for all}\quad u,v\in H^1(\Gamma).$$
Then we have for $u\in\dom(\Delta_\beta)$
$$\langle S''_\omega(\phi_{\beta})u,v\rangle_{(H^1)'\times H^1}=(L^\beta u, v)_2=\langle\mathfrak{J}L^\beta u,v\rangle_{(H^1)'\times H^1}.$$
It is obvious that $S''_\omega(\phi_{\beta})$ is the extension of $\mathfrak{J}L^\beta.$
\begin{definition}\label{linear_inst}
The standing wave  $e^{i\omega t}\phi_{\beta}(x)$ is said to be spectrally  unstable if there exist $\lambda$ with $\re \lambda>0$ and $w\in\dom(\Delta_\beta)$ such that 
\begin{equation*}\label{eigen}
-iS''_\omega(\phi_{\beta}) w=\lambda\mathfrak{J}w.
\end{equation*}

%Otherwise, the standing wave is said to be %linearly stable $e^{i\omega t}\phi_{\beta}$. 
\end{definition}
\noindent The notion of spectral  instability is particularly important since frequently its  presence leads to nonlinear instability. 
\begin{remark}
Notice that spectral instability implies 
that  $0$ is  unstable solution to the linearized equation
$$
i\partial_{t} v(t)=S''_{\omega}(\phi_\beta)v(t)
$$
in the sense of Lyapunov.  
\end{remark}

The search of $\lambda$ with $\re \lambda>0$
might be simplified essentially 
identifying $L^2_\mathbb{C}(\Gamma)$ with $L^2_\mathbb{R}(\Gamma)\oplus L^2_\mathbb{R}(\Gamma)$. Namely, take  $u=u_1+iu_2,\,\, u_1, u_2\in L^2_{\mathbb{R}}(\Gamma)$. Then one gets 
\begin{align*}
&L^\beta u=L^\beta_1u_1+iL^\beta_2u_2,\quad \dom(L^\beta_j)=\dom(\Delta_\beta),\\
& L^\beta_1v=-\Delta_\beta v+\omega v-p(\phi_{\beta})^{p-1}v, \quad L^\beta_2v=-\Delta_\beta v+\omega v-(\phi_{\beta})^{p-1}v.
\end{align*}
Thus, the operator $L^\beta$  in  $L^2_\mathbb{R}(\Gamma)\oplus L^2_\mathbb{R}(\Gamma)$   is interpreted as 
\begin{equation}\label{diag}
L^\beta=\left(\begin{array}{cc}
L_1^\beta&0\\ 
0& L_2^\beta
\end{array}\right).\end{equation}
 Observe that the operator $J$ of multiplication by $-i$ in $L^2_\mathbb{R}(\Gamma)\oplus L^2_\mathbb{R}(\Gamma)$ acts as 
$$J=\left(\begin{array}{cc}
0&I\\ 
-I& 0
\end{array}\right).$$ 
Hence we arrive at the equivalence  for $w\in\dom(\Delta_\beta)$
\begin{equation}\label{equiv}
-iS''_\omega(\phi_{\beta}) w=\lambda\mathfrak{J}w\quad \Longleftrightarrow \quad \mathfrak{J}\left(\begin{array}{cc}
0&L_2^\beta\\ 
-L_1^\beta& 0
\end{array}\right)\left(\begin{array}{c}
w_1\\
w_2
\end{array}\right)=\lambda\mathfrak{J}\left(\begin{array}{c}
w_1\\
w_2
\end{array}\right), \, w=w_1+iw_2.
\end{equation}
Our idea is to apply the approach by \cite{Grill88} for proving that the operator $\left(\begin{array}{cc}
0&L_2^\beta\\ 
-L_1^\beta& 0
\end{array}\right)$ has a positive  eigenvalue. 
To do that, we first need to study spectral properties of $L^\beta_1$ and $L^\beta_2$.
\begin{proposition}\label{spec_graph'}
Let $\beta\in \mathbb{R}\setminus\{0\}$, $\omega >\tfrac{N^2}{\beta^2}$, and $\tilde\phi_\beta$ be defined by \eqref{phi_beta}. Then the following assertions hold.  

\noindent$(i)$\,
$L^\beta_2\geq 0$  for $\beta<0$, and $\ker(L^\beta_2)=\Span\{\phi_{\beta}\}$ for any $\beta$.

\noindent$(ii)$\, If $\omega \neq
\tfrac{N^2}{\beta^2}\tfrac{p+1}{p-1}$, then
$\ker(L^\beta_1)=\{0\}$, while, for 
 $\omega =
\tfrac{N^2}{\beta^2}\tfrac{p+1}{p-1}$,
 the kernel of
$L^\beta_1$ is given by
$\ker(L^\beta_1)=\Span\{\hat{\phi}_1,\ldots,\hat{\phi}_{N-1}
\}$, where
\begin{equation*}\label{Psi}
\hat{\phi}_j=(0,\ldots,0,\underset{
j}{\tilde\phi_{\beta}},\underset{
j+1}{-\tilde\phi_{\beta}},0,\ldots,0).
 \end{equation*}
$(iii)$\,  $\sigma_{\ess}(L^\beta_1)=\sigma_{\ess}(L^\beta_2)=[\omega,\infty).$
\end{proposition}
\begin{proof}
$(i),(ii)$\,
Positivity of $L^\beta_2$ for $\beta<0$ and item   $(ii)$ were proven in \cite[Proposition 3.24]{AnGo2018}.
The identity $\ker(L^\beta_2)=\Span\{\phi_{\beta}\}$ follows from the fact that the only decaying solution to 
$$-v''+\omega v-(\phi_\beta)^{p-1}v=0$$
on $\Gamma$ is $\phi_\beta$ (up to a constant).  
 Indeed, to show the equality $\ker(L_2^\beta)=\Span\{\phi_\beta\}$,
we note that any  $v=(v_{j})_{j=1}^{N} \in \ker(L_2^\beta)$ has the form $v=(c_j \tilde\phi_\beta)_{j=1}^N$. From the inclusion $v\in\dom(L_2^\beta)$ we get $c_1=\ldots=c_N.$ 

 $(iii)$\, Consider the self-adjoint operator 
  $-\Delta_\infty=\bigoplus\limits_{j=1}^Nh_\infty,$
where \begin{equation*}
\begin{split}
(h_{\infty}{v})(x)=-v''(x),\quad x> 0,\quad
\dom(h_\infty)=\Big\{{v}\in H^2(\mathbb{R}^+): v'(0)=0\Big\}.
\end{split}
\end{equation*}
 Therefore, $\sigma_{\ess}(-\Delta_\infty)=\sigma_{\ess}(h_\infty)=[0,\infty).$
Notice that  the operators $-\Delta_\beta$  and  $-\Delta_\infty$ are self-adjoint extensions of the symmetric operator 
 \begin{equation*}
\begin{split}
(-\Delta_0 v)(x)&=\left(-v_j''(x)\right)_{j=1}^N,\quad x> 0,\quad  v=(v_j)_{j=1}^N,\\
\dom(\Delta_0)&=\Big\{v\in H^2(\Gamma): v'_1(0)=\ldots=v'_N(0)=0,\,\,\sum\limits_{j=1}^N  v_j(0)=0\Big\}.
\end{split}
\end{equation*}
The operator $-\Delta_0$ has equal deficiency indices $n_{\pm}(-\Delta_0)=\dim\ker(-\Delta_0^*\mp i)=1$, therefore, by Krein's resolvent formula, the operator $(-\Delta_\beta-\lambda)^{-1}- (-\Delta_\infty-\lambda)^{-1},\,\,\lambda\in\rho(-\Delta_\beta)\cap\rho(-\Delta_\infty),$ is of rank one (see \cite[Appendix A, Theorem A.2]{AlbGes05}). Then, by Weyl's theorem \cite[Theorem XIII.14]{ReeSim78},  $\sigma_{\ess}(-\Delta_\beta+\omega)=\sigma_{\ess}(-\Delta_\infty+\omega)=[\omega,\infty)$. The operator of multiplication by $(\phi_{\beta})^{p-1}$ is relatively $(-\Delta_\beta+\omega)$-compact,  therefore  $\sigma_{\ess}(L^\beta_1)=\sigma_{\ess}(L^\beta_2)=[\omega,\infty)$ (see Corollary 2 of \cite[Theorem XIII.14]{ReeSim78}).
 \end{proof}
 Further we study the number of negative eigenvalues of $L^\beta_1$ and  $L^\beta_2$.
\begin{proposition}\label{Morse_index}
$(i)$\, Let $\beta<0$, then 
$n(L^\beta_2)=0$. Moreover, 

\,\,  $1)$ if $\omega\leq \frac{p+1}{p-1}\frac{N^2}{\beta^2}$, then $n(L^\beta_1)=1$;
 
\,\, $2)$ if $\omega> \frac{p+1}{p-1}\frac{N^2}{\beta^2}$, then $n(L^\beta_1)=N.$
 
 $(ii)$\, Let $\beta>0$, then 
$n(L^\beta_2)=N-1$. Moreover, 

\,\, $1)$ if $\omega< \frac{p+1}{p-1}\frac{N^2}{\beta^2}$, then $n(L^\beta_1)=2N-1$;
 
\,\, $2)$ if $\omega\geq \frac{p+1}{p-1}\frac{N^2}{\beta^2}$, then $n(L^\beta_1)=N.$
 
\end{proposition}
\begin{proof}
We use a generalization of the
Sturm theory for the star graph elaborated in \cite{KaiPel18, Kai19}. 
Firstly, we calculate $n(L^\beta_1)$.   

\noindent\textit{Step 1.}\,  Suppose that $\lambda<0$ is an eigenvalue of $L^\beta_1$ with an   eigenvector $v^\lambda=(v^\lambda_j)_{j=1}^N\in \dom(L^\beta_1)$:\,  $L^\beta_1v^\lambda=\lambda v^\lambda.$ Then, denoting $a=\frac{2}{(p-1)\sqrt{\omega}}\tanh ^{-1}\left(\frac{N}{\beta \sqrt{\omega}}\right)$, we have 
$$-(v^\lambda_{j})''+\omega v^\lambda_{j}- \frac{p(p+1) \omega}{2} \operatorname{sech}^{2}\left(\frac{(p-1) \sqrt{\omega}}{2} (x-a)\right)v^\lambda_{j}=\lambda v^\lambda_{j}, \quad x \in(0, \infty).$$
Thus, $v_j^\lambda=c_ju(x-a),$ where $u(x)$ is the solution (on the line) to
$$-u''+\omega u- \frac{p(p+1) \omega}{2} \operatorname{sech}^{2}\left(\frac{(p-1) \sqrt{\omega}}{2} x\right)u=\lambda u.$$

By \cite[Lemma 4.1]{Kai19},  $u(x_{0})$ is a $C^{1}$ function of $\lambda$  for any fixed $x_{0} \in \mathbb{R},$  and 
$$\lim\limits_{\lambda \rightarrow-\infty}\frac{u'\left(x_{0}\right)}{u\left(x_{0}\right)}=-\infty, \quad \lim _{x \rightarrow+\infty} u(x) e^{\sqrt{\omega-\lambda} x}=1. $$
The coefficients $c_j$ satisfy the system 
 \begin{equation}\label{equation}
c_{1} u'(-a)=\ldots=c_{N} u'(-a), \quad  
\sum\limits_{j=1}^{N} c_{j} u(-a)=\beta c_{N} u'(-a).
\end{equation}
The determinant of the matrix associated with the  system
$$ \left(\begin{array}{ccccccccc}
u'(-a) & -u'(-a) & 0 & \ldots & 0 \\
u'(-a) & 0 & -u'(-a) & \ldots & 0  \\
\vdots & \vdots & \vdots & \ddots & \vdots   \\
u'(-a) & 0 & 0 & \ldots & -u'(-a)   \\
u(-a) & u(-a) & u(-a) & \ldots &  u(-a)-\beta u'(-a)
\end{array}\right)$$
 is given by 
\begin{equation}\label{D}
D=\left(u'(-a)\right)^{N-1}\left[N u(-a)-\beta u'(-a)\right].
\end{equation}
Hence we conclude that $\lambda$ is an eigenvalue of  $L^\beta_1$ if,  and only if,  either $u'(-a)=0$ (and multiplicity is $N-1$ in this case) or  $N u(-a)-\beta u'(-a)=0$.

Notice that two terms cannot be zero simultaneously  since $u(-a)=u'(-a)=0$ implies $u\equiv 0.$

\noindent\textit{Step 2.} We analyze  negative zeroes of $D$. Notice that $D$ can be expressed as $D=-\beta (u(-a, \lambda))^{N} F(\lambda)^{N-1}(F(\lambda)-\frac{N}{\beta}),$ where
 $$F(\lambda):=\dfrac{u'(-a, \lambda)}{u(-a, \lambda)}: (-\infty, 0]\to \mathbb{R}.$$
By \cite[Lemma 4.1, Lemma 4.5, Remark 4.6]{Kai19},  we get:

$\bullet$\, $\lim\limits_{\lambda\to -\infty}F(\lambda)=-\infty$;

$\bullet$\, if $\beta>0$, then there exists a unique pole $\lambda_*\in (-\infty, 0]$ of $F(\lambda)$ and 
$$\lim\limits_{\lambda\to \lambda_*^-}F(\lambda)=+\infty,\quad \lim\limits_{\lambda\to \lambda_*^+}F(\lambda)=-\infty,$$ 
moreover, $F(\lambda)$ is increasing on $(-\infty, \lambda_*)$ and $(\lambda_*, 0];$

$\bullet$\, if $\beta<0$, then  $F(\lambda)$  increases and is of class $C^1$ on $(-\infty, 0)$.

 It is easily seen that $u(x, 0)=\phi_\omega'(x),$ 
where $\phi_\omega(x)$ is defined by  \eqref{phi_omega},
 and $$F(0)-\frac{N}{\beta}=\frac{\beta \omega}{2N}\left(\frac{(p+1)N^2}{\beta^2\omega}-(p-1)\right)-\frac{N}{\beta}=\frac{\beta \omega(p-1)}{2N}\left(\frac{N^2}{\beta^2\omega}-1\right). $$

Since $\omega>\frac{N^2}{\beta^2}$, we get that  $F(0)>\frac{N}{\beta}$ for $\beta <0$, and $F(0)<\frac{N}{\beta}$ for $\beta >0$. Then, using properties of $F(\lambda)$ listed above,  we conclude that there exists a unique $\tilde\lambda<0$ such that $F(\tilde\lambda)-\frac{N}{\beta}=0.$ 
Moreover, we have (see the graph of $F(\tilde\lambda)$ for the different values of $\beta$ and $\omega$ on Figures 1-4 below):

$\bullet$\, $\beta>0\quad\Rightarrow\quad\left\{\begin{array}{c}
F(0)>0 \quad \text{for}\quad \omega<\frac{p+1}{p-1}\frac{N^2}{\beta^2},\\
F(0)\leq 0 \quad \text{for}\quad \omega\geq\frac{p+1}{p-1}\frac{N^2}{\beta^2}.
\end{array}\right. $ 

$\bullet$\, $\beta<0\quad\Rightarrow\quad\left\{\begin{array}{c}
F(0)\leq 0 \quad \text{for}\quad \omega\leq \frac{p+1}{p-1}\frac{N^2}{\beta^2},\\
F(0)>0 \quad \text{for}\quad \omega>\frac{p+1}{p-1}\frac{N^2}{\beta^2}.
\end{array}\right. $ 

\vspace{0.3cm}

\begin{tikzpicture}[every axis/.style={width=5cm}]
    \begin{axis}[axis lines=middle,ylabel=$y$,xlabel=$\lambda$,enlargelimits,  ytick={},
     yticklabels={},  xtick={},
     xticklabels={}, xmin=-4,xmax=1,
     ymin=-2,ymax=3]
      \addplot[thick,domain=-1.9:0,samples=200]{-1/(x+2)+2};
      \addplot[thick,domain=-4:-2.1,samples=200]{-2/(x+2)-2};
            
  ];

  \draw[thick,dashed,gray] (axis cs:-2,-4) -- (axis cs:-2,4);
  \draw[gray] (axis cs:-4,2) -- (axis cs:1,2);
  \draw[gray, dashed] (axis cs:-2.5,2) -- (axis cs:-2.5,0);
   \draw(1,2.5) node {\scalebox{0.8}{$y=\frac{N}{\beta}$}};
   \draw(-2.5,-0.5) node {\scalebox{0.9}{$\tilde\lambda$}};
 %  \draw(-1.2,-0.5) node {\scalebox{0.9}{$\lambda_2$}};
   \draw(-1.8,0.3) node {\scalebox{0.9}{$\lambda_*$}};
   
     \filldraw[blue] (-3,0) circle(2pt);
   \filldraw[blue] (-1.5,0) circle(2pt);
  \filldraw[blue] (-2.5,2) circle(2pt);
    \end{axis}
  
    \begin{axis}[xshift = 3in, axis lines=middle,ylabel=$y$,xlabel=$\lambda$,enlargelimits,  ytick={},
     yticklabels={},  xtick={},
     xticklabels={}, xmin=-4,xmax=1,
     ymin=-2,ymax=3]
      \addplot[thick,domain=-1.9:0,samples=200]{-1/(x+2)};
      \addplot[thick,domain=-4:-2.1,samples=200]{-2/(x+2)-2};
            
  ];

  \draw[thick,dashed,gray] (axis cs:-2,-4) -- (axis cs:-2,4);
  \draw[gray] (axis cs:-4,2) -- (axis cs:1,2);
  \draw[gray, dashed] (axis cs:-2.5,2) -- (axis cs:-2.5,0);
   \draw(1,2.5) node {\scalebox{0.8}{$y=\frac{N}{\beta}$}};
   \draw(-2.5,-0.5) node {\scalebox{0.9}{$\tilde\lambda$}};
      \draw(-1.8,0.3) node {\scalebox{0.9}{$\lambda_*$}};
   
   \filldraw[blue] (-2.5,2) circle(2pt);
   \filldraw[blue] (-3,0) circle(2pt);
    \end{axis}

    \draw(2,-0.5) node {\scalebox{0.7}{Figure 1: Graph of $F(\lambda)$ for $\beta>0$}} ; 
   \draw(2,-1) node {\scalebox{0.7}{and $\omega<\frac{p+1}{p-1}\frac{N^2}{\beta^2}$}}; 
   
   \draw(10,-0.5) node {\scalebox{0.7}{Figure 2: Graph of $F(\lambda)$ for $\beta>0$}} ; 
   \draw(10,-1) node {\scalebox{0.7}{and $\omega\geq\frac{p+1}{p-1}\frac{N^2}{\beta^2}$}};
     \end{tikzpicture}

\vspace{0.5cm}

\begin{tikzpicture}[every axis/.style={width=5cm}]
    \begin{axis}[axis lines=middle,ylabel=$y$,xlabel=$\lambda$,enlargelimits,  ytick={},
     yticklabels={},  xtick={},
     xticklabels={}, xmin=-4,xmax=1,
     ymin=-2,ymax=3]
      \addplot[thick,domain=-1.9:0,samples=200]{-1/(x+2)+2};
       ];

  \draw[gray] (axis cs:-4,-2) -- (axis cs:1,-2);
  \draw[gray, dashed] (axis cs:-1.75,-2) -- (axis cs:-1.75,0);
   \draw(1,-1.5) node {\scalebox{0.8}{$y=\frac{N}{\beta}$}};
   \draw(-1.9,0.5) node {\scalebox{0.9}{$\tilde\lambda$}};

   \filldraw[blue] (-1.5, 0) circle(2pt);
   \filldraw[blue] (-1.75,-2) circle(2pt);
    \end{axis}
  
    \begin{axis}[xshift = 3in, axis lines=middle,ylabel=$y$,xlabel=$\lambda$,enlargelimits,  ytick={},
     yticklabels={},  xtick={},
     xticklabels={}, xmin=-4,xmax=1,
     ymin=-3,ymax=2]
      \addplot[thick,domain=-1.9:0,samples=200]{-1/(x+2)-0.5};
       ];

  \draw[gray] (axis cs:-4,-2) -- (axis cs:1,-2);
  \draw[gray, dashed] (axis cs:-1.33,-2) -- (axis cs:-1.3,0);
   \draw(1,-1.5) node {\scalebox{0.8}{$y=\frac{N}{\beta}$}};
   \draw(-1.33,0.5) node {\scalebox{0.8}{$\tilde\lambda$}};
   
   \filldraw[blue] (-1.33,-2) circle(2pt);
    \end{axis}

    \draw(2,-0.5) node {\scalebox{0.7}{Figure 3: Graph of $F(\lambda)$ for $\beta<0$}} ; 
   \draw(2,-1) node {\scalebox{0.7}{and $\omega>\frac{p+1}{p-1}\frac{N^2}{\beta^2}$}}; 
   
   \draw(10,-0.5) node {\scalebox{0.7}{Figure 4: Graph of $F(\lambda)$ for $\beta<0$}} ; 
   \draw(10,-1) node {\scalebox{0.7}{and $\omega\leq\frac{p+1}{p-1}\frac{N^2}{\beta^2}$}};
     \end{tikzpicture}

Again, using properties of $F(\lambda),$ we conclude that:

$\bullet$\, $\beta<0\quad\Rightarrow$\quad there  exists a unique negative zero of  $F(\lambda)$ for $\omega>\frac{p+1}{p-1}\frac{N^2}{\beta^2}$, and for $\omega\leq \frac{p+1}{p-1} \frac{N^2}{\beta^2}$ we do not have negative zeroes of $F(\lambda)$.

$\bullet$\,  $\beta>0\quad\Rightarrow$\quad there  exists  a unique negative zero of  $F(\lambda)$ for $\omega\geq \frac{p+1}{p-1}\frac{N^2}{\beta^2}$, and for $\omega<\frac{p+1}{p-1}\frac{N^2}{\beta^2}$ we have two negative zeroes  of $F(\lambda).$

Finally, summarizing and noticing that negative zeros of $F(\lambda)$ are zeros of the determinant $D$ of multiplicity $N-1$, we conclude:

$\bullet$\, $\beta<0\quad\Rightarrow$\quad for $\omega\leq \frac{p+1}{p-1}\frac{N^2}{\beta^2}$   we have $n(L^\beta_1)=1$,  and for  $\omega>\frac{p+1}{p-1}\frac{N^2}{\beta^2}$ we have $n(L^\beta_1)=N$.  
 
$\bullet$\,  $\beta>0\quad\Rightarrow$\quad for $\omega< \frac{p+1}{p-1}\frac{N^2}{\beta^2}$   we have $n(L^\beta_1)=2N-1$,  and for  $\omega\geq\frac{p+1}{p-1}\frac{N^2}{\beta^2}$ we have $n(L^\beta_1)=N$.  

\noindent\textit{
Step 3.} We calculate $n(L^\beta_2)$ for $\beta>0$. Analogously to the previous case, the  number of negative eigenvalues of $L^\beta_2$ is determined by the number of zeroes of determinant \eqref{D}, where $u(x)$ is the solution to  
$$-u''+\omega u- \frac{(p+1) \omega}{2} \operatorname{sech}^{2}\left(\frac{(p-1) \sqrt{\omega}}{2} x\right)u=\lambda u.$$
In this case, the function $F(\lambda)$  increases and is of class  $C^1$ on $(-\infty, 0]$, and $\lim\limits_{\lambda\to -\infty}F(\lambda)=-\infty$.  
 The absence of the poles  follows from the Sturm oscillation  theorem since  $u(x, 0)=\phi_{\omega}(x)$ defined by \eqref{phi_omega}  (which is nonzero function).  The rest of the properties follow analogously (see \cite[Lemma 4.5]{Kai19}).
Observing that $F(0)=\frac{N}{\beta}$, we get the existence of  a unique negative zero of $F(\lambda),$ therefore,  $n(L^\beta_2)=N-1.$
\end{proof}
To prove orbital instability part in Theorem \ref{main_instub}, we apply the following abstract result by Henry, Perez, and Wreszinski \cite{HenPer82}. 
\begin{theorem}\label{Th_Henry}
 Let $X$ be a Banach space and let  $U \subset X$ be an open set containing $0$. Suppose that  $T: U \rightarrow X$ is such that  $T(0)=0$ and, for some $p>1$ and  $L\in B(X)$ with spectral radius $r(L)>1$, the relation holds 
$$
\|T(x)-L x\|_X=O\left(\|x\|^{p}_X\right) \quad \text { as } \quad x \rightarrow 0.
$$
Then $0$ is unstable as a fixed point of $T$: that is,  there is $\eta_0>0$ such that for all $\varepsilon>0$, there are $N_0\in\mathbb{N}$ and $x_0$ with $\|x_0\|_X\leq \varepsilon$  such that $\left\|T^{N_0}(x_0)\right\|_{X} \geq \eta_0.$
Moreover, if $\gamma$ is a $C^1$ curve of fixed points of $T$ such that $0\in \gamma$, then $\gamma$ is unstable.
\end{theorem}
The proof of the above theorem might be found in \cite[Section 2]{HenPer82}.
Actually on practice the corollary below appears to be useful.
\begin{corollary}\label{tilde_T}
Let $X$ be a Banach space,  $\phi\in X$, and let $\tilde U \subset X$ be  an open set containing $\phi$. Suppose that $\tilde T: \tilde U \rightarrow X$ is $C^2$ mapping satisfying $\tilde T(\phi)=\phi.$
  If there is an element $\mu \in \sigma\left(\tilde T^{\prime}(\phi)\right)$ with $|\mu|>1,$ then $\phi$ is an unstable fixed point of $\tilde T$. Moreover,   if $\gamma$ is a $C^1$ curve of fixed points of $\tilde T$ such that $\phi \in \gamma$, then $\gamma$ is unstable.
\end{corollary}
\begin{proof}
The proof was established in \cite[Corollary 3.1]{AngNat16}. We repeat it for convenience of the reader.

Let $x \in U:=\{y-\phi: y \in \tilde U\},$ and   $T(x):= \tilde T(x+\phi)-\phi$.  Then, we have $T(0)=\tilde T(\phi)-\phi=0$,  and $1<|\mu| \leq r(\tilde T'(\phi))$.
By the Taylor formula, $T(x)=T(0)+T^{\prime}(0) x+O\left(\|x\|_{X}^{2}\right)=\tilde T'(\phi)x+O\left(\|x\|_{X}^{2}\right)$ as  $x\to 0$. Therefore, by Theorem \ref{Th_Henry} with $L=\tilde T'(\phi)$, there is  $\eta_{0}>0$ such that for all $\varepsilon>0$ there exist  $N_0\in\mathbb{N}$ and $y_0$ with $\|y_0-\phi\|_X\leq \varepsilon$ such that $\left\|\tilde T^{N_0}(y_0)-\phi\right\|_{X} \geq \eta_{0} .$ To conclude, one observes that, by the last assertion of Theorem \ref{Th_Henry}, we get the existence of $\eta_0>0$ such that for all $\varepsilon>0$ there exist  $N_0\in\mathbb{N}$ and $y_0$ with $\|y_0-\phi\|_X\leq \varepsilon$ such that $\dist(\tilde T^{N_0}(y_0), \gamma)\geq \eta_0.$
\end{proof}
 \begin{proof}[Proof of Theorem \ref{main_instub}]
\textit{ Step 1.} Firstly, we comment on the proof of spectral instability. Observe that, by \eqref{equiv}, it is sufficient to prove that there exist $\lambda>0$ and $w_1, w_2\in \dom(\Delta_\beta)$ such that 
\begin{equation}\label{lambda_L}
\left(\begin{array}{cc}
0&L_2^\beta\\ 
-L_1^\beta& 0
\end{array}\right)\left(\begin{array}{c}
w_1\\
w_2
\end{array}\right)=\lambda\left(\begin{array}{c}
w_1\\
w_2
\end{array}\right).
\end{equation}
   \cite[Theorem 1.2]{Grill88} states that the number $I(L_1^\beta, L_2^\beta)$ of positive  $\lambda$ satisfying \eqref{lambda_L} is estimated by 
 \begin{equation}\label{I_estim}
 n(PL_1^\beta)-n(P(L_2^\beta)^{-1})\leq I(L_1^\beta, L_2^\beta),\end{equation}
 where $P$ is the orthogonal projection on $\ker(L^\beta_2)^\perp=\{v\in L^2(\Gamma):\, (v, \phi_{\beta})_2=0\}$ (here we assume that $L^2(\Gamma)$ is endowed with the usual complex scalar product).

By  spectral properties of $L^\beta_1, L^\beta_2$ proved in       Proposition  \ref{spec_graph'} and \ref{Morse_index}, we conclude the existence of such $\lambda>0$. 
  In particular, for $\beta$ and $\omega$ from the statement  of  Theorem \ref{main_instub} one gets $I(L_1^\beta, L_2^\beta)\geq N-1$.

\noindent\textit{Step 2.} To show orbital instability  for $p>2$ we use  Corollary \ref{tilde_T}. Let $u(t,x)$ be  a solution of \eqref{NLS}. Set $u(t,x)=e^{i\omega t}v(t,x)$, then $v(t,x)$ satisfies
\begin{equation}\label{NLS_modif}
i\partial_tv(t,x)=-\Delta_\beta v(t,x)+\omega v(t,x)-|v(t,x)|^{p-1}v(t,x)=S_\omega'(v(t,x)).
\end{equation}
Define $\tilde T: H^1(\Gamma)\rightarrow H^1(\Gamma)$ by $\tilde T(y)=v_y(\frac{2\pi}{\omega},x),$ where $v_y(t,x)$ is a solution of \eqref{NLS_modif} with initial value $v_y(0,x)=y(x)\in H^1(\Gamma)$. Notice that $\tilde T(\phi_{\beta})=\phi_{\beta}$ (since $\phi_{\beta}$ is the equilibrium point for \eqref{NLS_modif}). By Proposition \ref{PP1}-$(iii)$,  $\tilde T$ is $C^2$ mapping for $p>2.$ It is easily seen that $\tilde T'(\phi_{\beta})z=v_z(\frac{2\pi}{\omega},x)$, where $v_z(t,x)$ is the solution of the initial-value problem
\begin{equation*}
\partial_tv(t,x)=-iS''_\omega(\phi_{\beta})v(t,x), \quad v(0,x)=z(x).
\end{equation*}
By \textit{Step 1}, we conclude that  there exists $\lambda>0$ such that $\mu=e^\lambda\in \sigma(\tilde T'(\phi_{\beta}))$ (hence $|\mu|>1$). Thus, by Corollary \ref{tilde_T}, $\phi_\beta$ is nonlinearly  unstable, that is, there is  $\eta_{0}>0$ such that for all $\varepsilon>0$ there exist  $N_0$ and $y_0$ with $\|y_0-\phi_{\beta}\|_{H^1}\leq \varepsilon$ such that $\left\|\tilde T^{N_0}(y_0)-\phi_{\beta}\right\|_{H^1} \geq \eta_{0} ,$ or 
$$\|v_{y_0}(\tfrac{2\pi}{\omega} N_0)-\phi_{\beta}\|_{H^1}=\|u_{y_0}(\tfrac{2\pi}{\omega} N_0)-\phi_{\beta}\|_{H^1}\geq \eta_{0},$$
where $u_{y_0}(t,x)$ is the solution to \eqref{NLS} with initial value $y_0.$

Finally, observe that $\gamma=\{e^{i\theta}\phi_{\beta}: \theta\in [0, 2\pi)\}$ is $C^1$ curve of fixed points of $\tilde T$. Then $\gamma$ is nonlinearly unstable as well, that is, there is  $\eta_{0}>0$ such that for all $\varepsilon>0$  there exist  $N_0$ and $y_0$ with $\|y_0-\phi_{\beta}\|_{H^1}\leq \varepsilon$ such that 
\begin{align*}&\inf\limits_{\theta\in [0,2\pi)}\left\|\tilde T^{N_0}(y_0)-e^{i\theta}\phi_{\beta}\right\|_{H^1} =\inf\limits_{\theta\in [0,2\pi)}\|v_{y_0}(\tfrac{2\pi}{\omega} N_0)-e^{i\theta}\phi_{\beta}\|_{H^1}\\&=\inf\limits_{\theta\in [0,2\pi)}\|u_{y_0}(\tfrac{2\pi}{\omega} N_0)-e^{i\theta}\phi_{\beta}\|_{H^1}\geq \eta_{0}.
\end{align*} 
Thus, $e^{i\omega t}\phi_{\beta}$ is orbitally unstable. 
 \end{proof}
 \begin{remark}\label{remark_I}
 The idea of the proof of  \cite[Theorem 1.2]{Grill88}  is to consider the   solutions of the  problem 
$$\Big(PL^\beta_1+\lambda^2P(L^\beta_2)^{-1}\Big)v=0,\quad v\in \ker(L^\beta_2)^\perp,$$
instead of the ones of \eqref{lambda_L}.
The above equation appears naturally noticing that \eqref{lambda_L} is equivalent to $\left\{\begin{array}{c}
L_2^\beta w_2=\lambda w_1,\\
L_1^\beta w_1=-\lambda w_2.
\end{array}
\right.$
Since $\ran(L_2^\beta)\perp\ker(L_2^\beta),$ we get $w_1\in\ker(L_2^\beta)^{\perp}$. Hence $w_2=\lambda(L^\beta_2)^{-1}w_1$ and $L_1^\beta w_1+\lambda^2(L^\beta_2)^{-1}w_1=0.$
The projection $P$ serves to fit the problem into the Hilbert space $\ker(L^\beta_2)^\perp$.
 \end{remark}
  \iffalse
$\bullet$\,  $\sigma(L^\beta)=\sigma(L^\beta_1)\cup\sigma(L^\beta_2),$

$\bullet$\, $n(L^\beta)=n(L^\beta_1)+n(L^\beta_2),$

$\bullet$\,$\sigma_{\ess}(L^\beta)=\sigma_{\ess}(L^\beta_1)\cup\sigma_{\ess}(L^\beta_2)$.  
\fi
\begin{remark}
Consider the NLS-$\delta$ equation on the line:
 \begin{equation}\label{NLS_delta}
 i\partial_{t}u(t,x)=-\partial_x^2u(t,x)-\gamma\delta(x)u(t,x)-\left|u(t,x)\right|^{p-1}u(t,x),
 \end{equation}
 where $(t,x)\in \mathbb{R}\times\mathbb{R},\,\, \gamma\in \mathbb{R}\setminus\{0\}.$
 Equation \eqref{NLS_delta} has the standing wave solution $u(t,x)=e^{i\omega t}\phi_\gamma(x)$, where \newline
 $\phi_\gamma(x)=\left\{\frac{(p+1)\omega}{2}\mathrm{sech}^2\left(\frac{(p-1)\sqrt{\omega}}{2}|x|+\mathrm{arctanh}(\frac{\gamma}{2\sqrt{\omega}})\right)\right\}^{\frac{1}{p-1}}$.

Extensive study of the orbital stability of $e^{i\omega t}\phi_\gamma(x)$ was made in \cite{RJJ, FukOht08, CozFuk08}. The only case which was left open was for $3<p<5, \gamma<0,$ and $\omega=\omega_2>\frac{\gamma^2}{4}$, where $\partial_\omega\|\phi_\gamma\|^2_2\bigl|_{\omega=\omega_2}=0.$ 
Using the arguments from the proof of Theorem \ref{main_instub}, one may prove that  $e^{i\omega t}\phi_\gamma(x)$ is orbitally unstable in that case. 
Indeed, in \cite[Lemma 12]{CozFuk08} it was shown that $n(L^\gamma_1)=2$. Observing that $n(L^\gamma_2)=0$ and $\ker(L^\gamma_2)=\Span\{\phi_\gamma\}$ (the operators $L^\gamma_1, L^\gamma_2$  are the counterparts of $L^\beta_1, L^\beta_2$), by \eqref{I_estim}, we get $I(L_1^\gamma, L_2^\gamma)\geq 1$. Using the fact that   the mapping data-solution 
$u_0\in H^1(\mathbb{R})\mapsto u(t)\in C\left([0, T_0], H^1(\mathbb{R})\right)$  for \eqref{NLS_delta} is at least of class $C^2$ for $p>2$, we get the result.
\end{remark}
\subsection{Stability analysis for symmetric profile $\phi_{\beta}$  via Grillakis-Shatah-Strauss (GSS) approach in the case $\beta>0$}
The case of $\beta<0$ was studied in the framework of Grillakis-Shatah-Strauss approach
  \cite{AnGo2018}.
 In this subsection we complete the results of the previous subsection having repulsive $\delta'_s$ coupling. 
\begin{theorem}\label{GSS}
Let $\beta>0, \omega>\frac{N^2}{\beta^2}, \omega\neq \frac{p+1}{p-1}\frac{N^2}{\beta^2}, $ and $N\geq 2$, then 

%$(i)$\, for $1<p\leq 2, \omega<\frac{p+1}{p-1}\frac{N^2}{\beta^2}$ and even $N$ the standing wave $e^{i\omega t}\phi_\beta$ is orbitally unstable;

$(i)$\,  for $1<p\leq 3$ the standing wave $e^{i\omega t}\phi_\beta$ is orbitally stable in $H^1_{\eq}(\Gamma)$; 

$(ii)$\, for $3 <p< 5$ there exists $\omega^*>\frac{N^2}{\beta^2}$ such that  the standing wave $e^{i\omega t}\phi_\beta$ is orbitally unstable  in $H^1(\Gamma)$ as $\omega<\omega^*$ and orbitally stable in $H^1_{\eq}(\Gamma)$ as $\omega>\omega^*$;

$(iii)$\, for $p\geq 5$ the standing wave  $e^{i\omega t}\phi_\beta$ is orbitally unstable  in $H^1(\Gamma)$.
\end{theorem}
Let $\mc R: H^1(\Gamma)\to (H^1(\Gamma))'$ be the Riesz isomorphism. 
 Principal role in GSS approach  is  played  by the spectral properties of the  operator $\mc R^{-1}S''_\omega(\phi_{\beta}): H^1(\Gamma)\to H^1(\Gamma)$. We denote $\mc L^\beta:=\mc R^{-1}S''_\omega(\phi_{\beta})$.
  Since $S''_\omega(\phi_{\beta}):H^1(\Gamma)\to (H^1(\Gamma))'$ is bounded, the operator $ \mc L^\beta: H^1(\Gamma)\to H^1(\Gamma)$ is bounded and self-adjoint:
$$(\mc L^\beta u,v)_{H^1}=\langle S''_\omega(\phi_{\beta}) u,v\rangle _{(H^1)'\times H^1}=\langle S''_\omega(\phi_{\beta})v, u \rangle_{(H^1)'\times  H^1 }=(u, \mc L^\beta v)_{H^1}, \quad u,v\in H^1(\Gamma).$$
  Above we have used symmetry of  $S''_\omega(\phi_{\beta})$ and  the fact that $H^1(\Gamma)$ is the  real Hilbert space. 
  
 One of the crucial assumptions of GSS theory is the particular set of spectral properties of $\mc L^\beta$ (see \cite[Assumption 3]{GrilSha87}).  
  The following proposition links spectral properties of the operator $\mc L^\beta$ and operator $L^\beta$ associated with the bilinear form $\langle S''_\omega(\phi_{\beta}) u,v\rangle_{(H^1)'\times H^1}$ in $L^2(\Gamma)$ (for the proof see \cite[Lemma 4.5]{Stu08}). 
\begin{proposition}\label{spec_eq}
Let operator $L^\beta$ be defined by \eqref{L_b}, then 
\begin{align*}&\ker(\mc L^\beta)=\ker(L^\beta), \quad n(\mc L^\beta)=n(L^\beta),\,\text{ and}\\
&\inf \sigma_{\ess}(L^\beta)>0\,\, \Longrightarrow \,\,\inf \sigma_{\ess}(\mc L^\beta)>0.
\end{align*}
\end{proposition}
By \eqref{diag}, we conclude that to characterize  spectral properties of $\mc L^\beta$ it is sufficient to study spectral properties of $L^\beta_1, L^\beta_2.$
\begin{lemma}
Let $\beta>0, \omega>\frac{N^2}{\beta^2}$, then  $n(L^\beta_1)=1$ and  $L^\beta_2\geq 0$   in $L^2_{\eq}(\Gamma)=\{v\in L^2(\Gamma):\, v_1(x)=\ldots =v_N(x)\}.$ 
\end{lemma}
\begin{proof}
To prove $n(L^\beta_1)=1$, notice that in $L^2_{\eq}(\Gamma)$ conditions \eqref{equation} turn into $Nu(-a)=\beta u'(-a)$, i.e.  $\lambda<0$ is an eigenvalue of $L^\beta_1$ iff  $F(\lambda)=\frac{N}{\beta}.$
 By the properties of $F(\lambda)$ for $\beta>0$ and since $F(0)<\frac{N}{\beta}$, we conclude that there exists a unique $\lambda_{\eq}<0$ such that $F(\lambda_{\eq})=\frac{N}{\beta}.$
 
 Identity $n(L^\beta_2)=0$  follows by the same reasoning and by the fact that $F(\lambda)$ is the increasing function of class $C^1$  on $(-\infty, 0]$ (see \textit{Step 3} in the proof of Proposition \ref{Morse_index}). Finally, by $\sigma_{\ess}(L^\beta_2)=[\omega, \infty),$ we get positivity of $L^\beta_2.$
\end{proof}
By Proposition \ref{spec_eq}, we get $n(\mc L^\beta)=1$, while, by Proposition \ref{spec_graph'}, $\ker(\mc L^\beta)=\Span\{i\phi_\beta\}$ for $\omega\neq \frac{p+1}{p-1}\frac{N^2}{\beta^2}$, and $\sigma_{\ess}(\mc L^\beta)=[\omega, \infty)$. Thus,   \textit{Assumption 3} in \cite{GrilSha87} holds for $\omega\neq \frac{p+1}{p-1}\frac{N^2}{\beta^2}$, and from \cite[Theorem 3]{GrilSha87},  we get:
\begin{theorem}\label{S_sinal}
The standing wave solution  $e^{i\omega t}\phi_\beta$ is orbitally stable  in $H^1_{\eq}(\Gamma)$ if, and only if, $\partial_\omega^2 S_\omega(\phi_\beta)>0$. 
\end{theorem}
Below we study the sign of $\partial_\omega^2 S_\omega(\phi_\beta)$.
\begin{proposition}\label{sinal}
Let $\omega>\frac{N^{2}}{\beta^{2}}, \beta>0$, and $J(\omega)=\partial_\omega^2 S_\omega(\phi_\beta)$.

$(i)$\, If $1<p \leq 3$, then $J(\omega)>0$.

$(ii)$\, If $3<p<5$, then there exists $\omega^{*}>\frac{N^{2}}{\beta^{2}}$ such that $J\left(\omega^{*}\right)=0$, and $J(\omega)<0$ for $\omega< \omega^*$, while $J(\omega)>0$ for $\omega>\omega^{*}.$

$(iii)$\, If $p\geq 5$, then $J(\omega)<0$.
\end{proposition}
\begin{proof}
Notice that due to $S'_\omega(\phi_\beta)=0,$ we get $J(\omega)=\partial_\omega\|\phi_\beta\|_2^2$.
 Recall that $\phi_\beta=(\tilde\phi_\beta)_{j=1}^{N}$  is defined by \eqref{phi_beta}. Then, via change of variables,
we get
$$
\int_{0}^{\infty}\tilde\phi^2_\beta(x) d x=\left(\frac{p+1}{2}\right)^{\frac{2}{p-1}} \frac{2 \omega^{\frac{2}{p-1}-\frac{1}{2}}}{p-1} \int_{\frac{N}{-\beta \sqrt{\omega}}}^{1}\left(1-t^{2}\right)^{\frac{2}{p-1}-1} d t .
$$
The last  equality yields 
$$
J(\omega)=C \omega^{\frac{7-3 p}{2(p-1)}} J_{1}(\omega), \quad C=\frac{N}{p-1}\left(\frac{p+1}{2}\right)^{\frac{2}{p-1}},
$$
where
$$
J_{1}(\omega)=\frac{5-p}{p-1} \int_{\frac{N}{-\beta \sqrt{\omega}}}^{1}\left(1-t^{2}\right)^{\frac{3-p}{p-1}} d t-\frac{N}{\beta \sqrt{\omega}}\left(1-\frac{N^{2}}{\beta^{2} \omega}\right)^{\frac{3-p}{p-1}} .
$$
Thus,
$$
J_{1}^{\prime}(\omega)=-\frac{N}{\beta \omega^{3 / 2}} \frac{3-p}{p-1}\left[\left(1-\frac{N^{2}}{\beta^{2} \omega}\right)^{\frac{3-p}{p-1}}+\frac{N^{2}}{\beta^{2} \omega}\left(1-\frac{N^{2}}{\beta^{2} \omega}\right)^{-\frac{2(p-2)}{p-1}}\right].
$$
It is easily seen that $J'_1(\omega)< 0$ as $1<p<3$, and $J'_1(\omega)> 0$ as $p>3$.
Observe that for $p<5$
\begin{align*}
\lim\limits_{\omega\to\infty}J_1(\omega)=\frac{5-p}{p-1}\int_0^1(1-t^2)^{\frac{3-p}{p-1}}dt=:a_0>0,\quad 
\lim\limits_{\omega\to\frac{N^2}{\beta^2}}J_1(\omega)=\left\{\begin{array}{c}
2a_0,\, p\in(1, 3],\\
-\infty, p\in (3, 5).
\end{array}\right.
\end{align*}
 By the above properties of $J_1(\omega),$ we conclude $J_1(\omega)>0$ for $1<p\leq 3$, and there exists $\omega^{*}>\frac{N^{2}}{\beta^{2}}$ such that $J_1(\omega^*)=0$, and 
 $J_1(\omega)<0$ for $\omega< \omega^*$, while $J_1(\omega)>0$ for $\omega>\omega^{*}.$
 Observing that $J_1(\omega)<0$ for $p\geq 5$, we conclude the proof.
\end{proof}
Now Theorem \ref{GSS} follows from  Proposition \ref{sinal} and Theorem \ref{S_sinal}. 
\begin{remark}\label{GSS_N_2}
Let $N$ be even and $\beta<0$. Using  the fact that $\phi_{\frac{N}{2}}$ is the minimizer of problem \eqref{min_N_2} for $\omega>\frac{p+1}{p-1}\frac{N^2}{\beta^2}$, we get from \cite[Proposition 6.11]{AdaNoj13} (using GSS approach) the following stability/instability result:

$(i)$\,  If $1<p\leq 5$, then for any $\omega>\frac{p+1}{p-1}\frac{N^2}{\beta^2}$, the standing wave $e^{i\omega t}\phi_{\frac{N}{2}}$ is orbitally stable in $H^1_{\frac{N}{2}}(\Gamma)$.

$(ii)$\, If $p>5$, then

\,\,$1)$ there exists $\omega_{1}$ such that $e^{i\omega t}\phi_{\frac{N}{2}}$ is orbitally stable in $H^1_{\frac{N}{2}}(\Gamma)$ for  $\frac{p+1}{p-1}\frac{N^2}{\beta^2}<\omega<\omega_{1}$;
 
\,\,$2)$ there exists $\omega_{2} \geq \omega_{1}$ such that $e^{i\omega t}\phi_{\frac{N}{2}}$ is orbitally unstable in $H^1_{\frac{N}{2}}(\Gamma)$ (and consequently in $H^1(\Gamma)$) for $\omega>\omega_{2}$. 
\end{remark}
\subsection{Strong instability results for $\phi_\beta$ and $\phi_{\frac{N}{2}}$}
In this subsection we complete the previous results introducing  particular type of the orbital instability by blow up. The formal definition is the following.
\begin{definition}
Standing wave $e^{i\omega t}\phi(x)$ is strongly unstable if for any $\varepsilon>0$ there exists $u_{0} \in H^1(\Gamma)$ such that $\left\|u_{0}-\phi(x)\right\|_{H^{1}}<\varepsilon$ and the solution $u(t)$ of \eqref{NLS} with $u(0)=u_{0}$ blows up in finite time.
\end{definition}
\noindent The strong instability results  are the next two theorems.
\begin{theorem}\label{beta>0}
Let  $\beta>0$, $\omega>\frac{N^2}{\beta^2}$, $N\geq 2,$ and $p\geq 5$, then  
the standing wave $e^{i\omega t}\phi_\beta(x)$ is  strongly unstable in $H^1(\Gamma)$.
 \end{theorem}
  
\begin{theorem}\label{beta<0}
Let $\beta<0,\, p>5,\omega>\frac{N^2}{\beta^2}$. Then the following assertions hold.

$(i)$\, Let $N\geq 2$. Assume that  $\hat{\xi}(p)\in(0,1)$ is a unique solution to
$$\frac{(p-5)N}{2}\int_\xi^1(1-s^2)^{\frac{2}{p-1}}ds=\xi(1-\xi^2)^{\frac{2}{p-1}},\quad (0<\xi<1),$$ and define $\omega_3=\omega_3(p,\beta)=\frac{N^2}{\beta^2\hat{\xi} ^2(p)}$. 
Then the standing wave solution $e^{i\omega t}\phi_\beta(x)$ is strongly unstable in $H^1(\Gamma)$ for all $\omega\in [\omega_3, \infty)$. 

$(ii)$\, Let $N$ be even, then there exists $\omega_4>\frac{p+1}{p-1}\frac{N^2}{\beta^2}$ such that $e^{i\omega t}\phi_{\frac{N}{2}}(x)$ is strongly unstable in $H^1(\Gamma)$ for all $\omega\in [\omega_4, \infty)$. 
\end{theorem}
\begin{remark}
Observe that  Theorem \ref{beta<0}-$(ii)$ completes orbital instability result   of Remark \ref{GSS_N_2}-$(ii),2)$.
\end{remark}
\noindent The proofs of the above theorems  are analogous to  the ones of \cite[Theorem 1.3 and Theorem 1.4]{NataMasa2020}.  The principal ingredients are  virial identity \eqref{well_17a} and Lemmas \ref{lemma_eq},  \ref{N/2}-$(ii).$
\begin{remark}
Analogously to \cite[Theorem 1.4]{ArdCel21} (i.e. using virial identity and variational characterization) one can prove for  the ground state solution $\phi_{\beta,V}$ of equation \eqref{V_beta}  the following result. Assume that   $p>5, \beta\in (\beta^*, -\frac{N}{\sqrt{\omega}}),  \omega>-\inf\sigma(-\Delta_\beta+V)$, and $V(x)$ satisfies assumptions \eqref{V_cond} and $xV'(x)\in L^1(\Gamma)+L^{\infty}(\Gamma)$, then  condition  $\partial_{\lambda}^{2} E\left(\phi_{\beta, V}^{\lambda}\right)|_{\lambda=1}<0$ implies orbital instability of the  standing wave solution  $e^{i \omega t} \phi_{\beta, V}(x)$ of \eqref{V_beta}  in $H^{1}(\Gamma)$.
 
\end{remark}

\subsection{Instability analysis for the asymmetric  profiles $\phi_k$}\label{assym}
The principal result of this subsection is the following theorem. 
\begin{theorem}\label{main_asym}
Let $p>1$, $\omega>\frac{p+1}{p-1}\frac{N^2}{\beta^2}$, and let  the profiles $\phi_k$ be defined in Theorem \ref{phi_k} for $k=1, \ldots, N-1$.

$(i)$\, If $\beta<0$, then $e^{i\omega t}\phi_k$ is spectrally  unstable for $k\geq 2$.

$(ii)$\,  If $\beta>0$, then $e^{i\omega t}\phi_k$ is spectrally unstable for $N-k\geq 4.$

If, additionally,  $p>2$, then in all the cases mentioned above we have orbital instability.

\end{theorem}
\noindent To prove the above theorem we repeat algorithm from the previous subsection, that is, firstly we linearize \eqref{NLS} at $\phi_k$. 
We denote the operators  analogous to $L^\beta_j$ by $L^{\phi_k}_j$, and  $\dom(L_j^{\phi_k})=\dom(\Delta_\beta)$.
 The proof of the proposition below repeats the one of Proposition \ref{spec_graph'}-$(i),(ii)$. 
\begin{proposition}\label{spec_L_assym}
Let $\beta\in \mathbb{R}\setminus\{0\}$ and $\omega >\frac{p+1}{p-1}\tfrac{N^2}{\beta^2}$, then the following assertions hold.  

\noindent$(i)$\,
$\ker(L^{\phi_k}_2)=\Span\{\phi_{k}\}$.

\noindent$(ii)$\,  $\sigma_{\ess}(L^{\phi_k}_j)=[\omega,\infty), \,\, j=1,2$.
\end{proposition}
\noindent Next we study the Morse index of $L^{\phi_k}_1$ and $L^{\phi_k}_2$.
  \begin{proposition}\label{Morse_index_assym} Suppose that $\omega>\frac{p+1}{p-1}\tfrac{N^2}{\beta^2}$ and $N\geq 2$. Then

$(i)$ for $\beta>0$ we have
$n(L^{\phi_k}_2)\leq N-1$ and $n(L^{\phi_k}_1)\geq 2N-k-3;$

$(ii)$\, for $\beta<0$ we have  $L^{\phi_k}_2\geq 0$  and 
$n(L^{\phi_k}_1)\geq k$.
\end{proposition}
\begin{proof}
Mainly we will use the ideas from the proof of Proposition \ref{Morse_index}.

\noindent\textit{Step 1.}    Suppose that $\lambda<0$ is an eigenvalue of $L^{\phi_k}_1$ with an   eigenvector $v^\lambda=(v^\lambda_j)\in \dom(L^{\phi_k}_1)$:\,  $L^{\phi_k}_1v^\lambda=\lambda v^\lambda.$ Then, denoting
 $$a_1=-\sign(\beta)\frac{2}{(p-1)\sqrt{\omega}}\tanh ^{-1}(t_1),\quad  a_N=-\sign(\beta)\frac{2}{(p-1)\sqrt{\omega}}\tanh ^{-1}(t_N),$$  and rearranging edges of $\Gamma$ (putting $k$ edges corresponding to $t_1$ at the beginning),  we get for $x>0$
\begin{align*}
&-(v^\lambda_{j})''+\omega v^\lambda_{j}- \frac{p(p+1) \omega}{2} \operatorname{sech}^{2}\left(\frac{(p-1) \sqrt{\omega}}{2} (x+a_1)\right)v^\lambda_{j}=\lambda v^\lambda_{j},\quad j=1,\ldots,k, \\
&-(v^\lambda_{j})''+\omega v^\lambda_{j}- \frac{p(p+1) \omega}{2} \operatorname{sech}^{2}\left(\frac{(p-1) \sqrt{\omega}}{2} (x+a_N)\right)v^\lambda_{j}=\lambda v^\lambda_{j},\quad j=k+1,\ldots,N.
\end{align*}
Thus, $v_j^\lambda=\left\{\begin{array}{c}
 c_ju(x+a_1),\, j=1,\ldots, k,\qquad\\
 c_j u(x+a_N),\, j=k+1,\ldots, N,  
\end{array}\right.$ where $u(x)$ is the solution  on the line to
$$-u''+\omega u- \frac{p(p+1) \omega}{2} \operatorname{sech}^{2}\left(\frac{(p-1) \sqrt{\omega}}{2} x\right)u=\lambda u.$$
The coefficients $c_j$ satisfy the system 
 \begin{align*}
&c_{1} u'(a_1)=\ldots=c_{k} u'(a_1)=c_{k+1}u'(a_N)=\ldots=c_N u'(a_N),\\
&\sum\limits_{j=1}^{k} c_{j} u(a_1)+\sum\limits_{j=k+1}^{N} c_{j} u(a_N)=\beta c_{N} u'(a_N).\end{align*}
The determinant of the matrix associated with the  system 
$$\left(\begin{array}{ccccccccc}
u'(a_1) & -u'(a_1) & 0 & \ldots & 0 & 0 & \ldots & 0 & 0 \\
u'(a_1) & 0 & -u'(a_1) & \ldots & 0 & 0 & \ldots & 0 & 0 \\
\vdots & \vdots & \vdots & \ddots & \vdots & \vdots & \ddots & \vdots & \vdots \\
u'(a_1) & 0 & 0 & \ldots & -u'(a_N) & 0 & \ldots & 0 & 0 \\
u'(a_1) & 0 & 0 & \ldots & 0 & -u'(a_N) & \ldots & 0 & 0 \\
\vdots & \vdots & \vdots & \ddots & \vdots & \vdots & \ddots & \vdots & \vdots \\
u'(a_1) & 0 & 0 & \ldots & 0 & 0 & \ldots & 0 & -u'(a_N) \\
u(a_1) & u(a_1) & u(a_1) & \ldots & u(a_1) & u(a_N) & \ldots & u(a_N) & u(a_N)-\beta u'(a_N)
\end{array}\right)$$
is given by 
\begin{equation}\label{D_asym}
D=(u'(a_1))^{k-1}\left(u'(a_N)\right)^{N-k-1}\Big[ku(a_1)u'(u_N)+(N-k)u(a_N)u'(a_1)-\beta u'(a_N)u'(a_1)\Big].
\end{equation}
Hence we conclude that $\lambda$ is an eigenvalue of  $L^{\phi_k}_1$ if,  and only if,  either $u'(a_1)=0$  (with multiplicity  $k-1$), or  $u'(a_N)=0$  (with multiplicity  $N-k-1$), or $ku(a_1)u'(a_N)+(N-k)u(a_N)u'(a_1)-\beta u'(a_N)u'(a_1)=0$.

We rewrite $D$ as 
\begin{equation}\label{D_as_1}
D=\Big(F_1(\lambda)u(a_1)\Big)^{k-1}\Big(F_N(\lambda)u(a_N)\Big)^{N-k-1}\Big[ku(a_1)u'(a_N)+(N-k)u(a_N)u'(a_1)-\beta u'(a_N)u'(a_1)\Big], 
\end{equation}
where 
$$F_1(\lambda):=\dfrac{u'(a_1, \lambda)}{u(a_1, \lambda)}:(-\infty, 0]\to \mathbb{R},\quad F_N(\lambda):=\dfrac{u'(a_N, \lambda)}{u(a_N, \lambda)}: (-\infty, 0]\to \mathbb{R}.$$
By \cite[Lemma 4.1, Lemma 4.5, Remark 4.6]{Kai19},  we get for $ j=1,N$:

$\bullet$\, $\lim\limits_{\lambda\to -\infty}F_j(\lambda)=-\infty$.

$\bullet$\, If $\beta>0$, then there exists a unique pole $\lambda_*^j\in (-\infty, 0]$ and 
$$\lim\limits_{\lambda\to {\lambda_*^j}-}F_j(\lambda)=+\infty,\quad \lim\limits_{\lambda\to \lambda_*^j +}F_j (\lambda)=-\infty,$$ 
moreover, $F_j(\lambda)$ is increasing on $(-\infty, \lambda_*^j)$ and $(\lambda_*^j, 0].$

$\bullet$\, If $\beta<0$, then  $F_j(\lambda)$  increases and is of class $C^1$ on $(-\infty, 0)$.
 
As in the previous case, observing that $u(x,0)=\phi'_\omega(x),$  we get 
\begin{equation*}
F_j(0)=\frac{\sqrt{\omega}}{2\sign(\beta)(p+1)t_j}\left(t_j^2-\frac{p-1}{p+1}\right).
\end{equation*}
Since $t_1<t_N$,   we conclude:

$\bullet$\, $\beta<0\quad\Rightarrow\quad\left\{\begin{array}{c}
F_1(0)>0,\\
F_N(0)<0.
\end{array}\right. $

$\bullet$\, $\beta>0\quad\Rightarrow\quad\left\{\begin{array}{c}
F_1(0)<0,\\
F_N(0)>0.
\end{array}\right. $

\noindent Using the properties of $F_j(\lambda),$ we obtain:

$\bullet$\, $\beta<0\quad\Rightarrow$\quad there  exists a unique $\lambda_1^1<0$ such that $F_1(\lambda_1^1)=0$, and $F_N(\lambda)$ does not have negative zeroes.  

$\bullet$\,  $\beta>0\quad\Rightarrow$\quad there  exists a unique $\lambda_1^1<0$ (with $\lambda_1^1<\lambda_*^1$) such that $F_1(\lambda_1^1)=0$, and $F_N(\lambda)$ has  two negative zeroes $\lambda_1^N, \lambda_2^N<0$ (with $\lambda_1^N<\lambda_*^N<\lambda_2^N$).

Now let $\beta<0$  and $\lambda^1_1$ be a unique negative zero of $F_1(\lambda)$. Consider the function 
$$\tilde F(\lambda)=\frac{k}{F_1(\lambda)}+\frac{N-k}{F_N(\lambda)}:(-\infty, \lambda_1^1)\rightarrow \mathbb{R}.$$
From the properties of $F_1(\lambda)$ and $F_N(\lambda),$ we conclude that $\tilde F(\lambda)$ is decreasing and 

$$\lim\limits_{\lambda\to -\infty}\tilde F(\lambda)=0-,\quad \lim\limits_{\lambda\to \lambda_1^1-}\tilde F (\lambda)=-\infty.$$
Therefore, there exists $\tilde\lambda<\lambda_1^1$ such that $\tilde F(\tilde \lambda)=\beta,$ and the last term in $D$ is zero at $\tilde\lambda$ (see the graph of the function $\tilde F(\lambda)$ on Figure 5).

\begin{tikzpicture}[every axis/.style={width=6cm}]
    \begin{axis}[axis lines=middle,ylabel=$y$,xlabel=$\lambda$,enlargelimits,  ytick={},
     yticklabels={},  xtick={},
     xticklabels={}, xmin=-4,xmax=1,
     ymin=-3,ymax=1.5]
      \addplot[thick,domain=-4:-1.01,samples=200]{1/(x+1)};
     
  ];

  \draw[thick,dashed,gray] (axis cs:-1,-4) -- (axis cs:-1,2);
  \draw[gray] (axis cs:-4,-2) -- (axis cs:4,-2);
  \draw[gray, dashed] (axis cs:-1.5,-2) -- (axis cs:-1.5,0);
   \draw(1,-1.5) node {\scalebox{0.8}{$y=\beta$}};

  \draw(-0.8,0.5) node {\scalebox{0.9}{$\lambda_1^1$}};
   \draw(-1.5,0.5) node {\scalebox{0.9}{$\tilde\lambda$}};
   
   \filldraw[blue] (-1.5,-2) circle(2pt);
 
    \end{axis}
  
    \begin{axis}[xshift = 3in, axis lines=middle,ylabel=$y$,xlabel=$\lambda$,enlargelimits,  ytick={},
     yticklabels={},  xtick={},
     xticklabels={}, xmin=-4,xmax=1,
     ymin=-3,ymax=3]
      \addplot[thick,domain=-4:-1.01,samples=200]{1/(x+1)};
        \addplot[thick,domain=-0.99:1,samples=200]{1/(2*x+2)};
      
      \addplot[thick,dashed,domain=-4:-2.1,samples=200]{1/(x+2)};
         \addplot[thick,dashed,domain=-1.9:1,samples=200]{1/(2*x+4)};      
  ];

  \draw[thick,dashed,gray] (axis cs:-2,-4) -- (axis cs:-2,4);
    \draw[thick,dashed,gray] (axis cs:-1,-4) -- (axis cs:-1,4);
 % \draw[gray] (axis cs:-4,-2) -- (axis cs:2,-2);

   \draw(-2.3,0.5) node {\scalebox{0.9}{$\lambda_1^1$}};
      \draw(-1.3,-0.5) node {\scalebox{0.9}{$\lambda_1^N$}};
   
    \end{axis}

    \draw(3,-1) node {\scalebox{0.7}{Figure 5: Graph of $\tilde F(\lambda)$  for $\beta<0$}} ; 
   \draw(3,-1.5) node {\scalebox{0.7}{on $(-\infty, \lambda_1^1)$}}; 
   
   \draw(11,-1) node {\scalebox{0.7}{Figure 6: The case of $\beta>0$ and $L_2^{\phi_k}$.}} ;
    
   \draw(11,-1.5) node {\scalebox{0.7}{Solid line is the graph of $1/F_N(\lambda),$ and }} ; 
   \draw(11,-2) node {\scalebox{0.7}{ the dashed line is the graph of  $1/F_1(\lambda)$}};
     \end{tikzpicture}

Summarizing the above, by \eqref{D_as_1}, we deduce:

$\bullet$\, $\beta<0\quad\Rightarrow$\quad $n(L^{\phi_k}_1)\geq k$.

$\bullet$\, $\beta>0\quad\Rightarrow$\quad
 $n(L^{\phi_k}_1)\geq 2N-k-3$.

\noindent\textit{Step 2.} We calculate $n(L^{\phi_k}_2)$  for $\beta>0$.
Analogously to the previous case, the  number of negative eigenvalues of $L^{\phi_k}_2$ is determined by the number of zeroes of the determinant \eqref{D_asym}, where $u(x)$ is the solution to  
$$-u''+\omega u- \frac{(p+1) \omega}{2} \operatorname{sech}^{2}\left(\frac{(p-1) \sqrt{\omega}}{2} x\right)u=\lambda u.$$
In this case, the function $F_j(\lambda), j=1,N,$  increases and is of class  $C^1$ on $(-\infty, 0]$ and $\lim\limits_{\lambda\to -\infty}F_j(\lambda)=-\infty$ (this also holds for $\beta<0$).  
 The absence of the poles  follows from the Sturm oscillation  theorem since  $u(x, 0)=\phi_{\omega}$. % where $\phi_\omega(x)=\left\{\tfrac{(p+1)\omega}{2}\mathrm{sech}^2\left(\tfrac{(p-1)\sqrt{\omega}}{2}x\right)\right\}^{\frac{1}{p-1}}$ (which is nonzero function). 
 Observe that $F_j(0)=\sqrt{\omega}t_j>0$, then there exists a unique $\lambda_1^j<0$  such that  $F_j(\lambda_1^j)=0.$ Without loss of generality we may assume that $\lambda_1^1<\lambda_1^N$. 
 Consider 
$$\tilde F(\lambda)=\frac{k}{F_1(\lambda)}+\frac{N-k}{F_N(\lambda)}:(-\infty, \lambda_1^1)\cup (\lambda_1^1, \lambda_1^N)\cup (\lambda_1^N, 0]\rightarrow \mathbb{R}.$$

The function is decreasing on each interval of its domain, moreover it satisfies (see the graphs  of the functions $1/F_1(\lambda)$ and  $1/F_N(\lambda)$ on Figure 6):
\begin{align*}
&\lim\limits_{\lambda\to -\infty}\tilde F(\lambda)=0-,\quad \lim\limits_{\lambda\to \lambda_1^1-}\tilde F (\lambda)=-\infty,\quad \lim\limits_{\lambda\to \lambda_1^1+}\tilde F (\lambda)=+\infty,\\
 & \lim\limits_{\lambda\to \lambda_1^N-}\tilde F (\lambda)=-\infty, \quad  \lim\limits_{\lambda\to \lambda_1^N+}\tilde F (\lambda)=+\infty, \quad \tilde F(0)=\frac{1}{\sqrt{\omega}}\left(\frac{k}{t_1}+\frac{N-k}{t_N}\right)=\beta.
\end{align*}  

From the above properties of $\tilde F(\lambda)$ it follows that there might exist $\tilde \lambda\in(\lambda_1^1, \lambda_1^N)$ such that $\tilde F(\tilde \lambda)=\beta$, and intervals $(-\infty, \lambda_1^1)$ and $ (\lambda_1^N, 0]$ do not contain zeroes of $\tilde F(\lambda)-\beta$.
Therefore, taking into account multiplicity of $\lambda_1^j$, we obtain $n(L^{\phi_k}_2)\leq N-1.$

For $\beta<0$ the identity $n(L^{\phi_k}_2)=0$ follows observing that  $F_j(0)=-\sqrt{\omega}t_j<0$ and $F_j(\lambda)$ is  increasing on $(-\infty, 0]$. In particular, there does not exist any negative zero of $\tilde F(\lambda)-\beta$. Indeed, $\tilde F(\lambda)$  is decreasing on $(-\infty, 0]$ and 
$$\lim\limits_{\lambda\to -\infty}\tilde F(\lambda)=0-, \quad \tilde F(0)=\frac{-1}{\sqrt{\omega}}\left(\frac{k}{t_1}+\frac{N-k}{t_N}\right)=\beta.$$
 Recalling   $\sigma_{\ess}(L^{\phi_k}_2)=[\omega, \infty)$, we get the positivity of $L^{\phi_k}_2$. 
\end{proof}
\begin{remark}
For $\beta<0$ and symmetric profile $\phi_\beta$,  positivity of the operator $L^\beta_2$ was proven in \cite[Proposition 3.24]{AnGo2018} using the  Jensen inequality for the function $f(x)=x^2$.
\end{remark}
\begin{proof}[ Proof of Theorem \ref{main_asym}]
As in the previous case, the principal ingredient of the proof is  to check the existence  of
 $\lambda>0$ and $w_1, w_2\in \dom(\Delta_\beta)$ such that 
\begin{equation}\label{lambda_L_asym}
\left(\begin{array}{cc}
0&L_2^{\phi_k}\\ 
-L_1^{\phi_k}& 0
\end{array}\right)\left(\begin{array}{c}
w_1\\
w_2
\end{array}\right)=\lambda\left(\begin{array}{c}
w_1\\
w_2
\end{array}\right).
\end{equation}
By  \cite[Theorem 1.2]{Grill88} and the  spectral properties of $L^{\phi_k}_1, L^{\phi_k}_2$ proved in       Proposition  \ref{spec_L_assym} and \ref{Morse_index_assym}, we conclude the existence of such $\lambda>0$   for $k$ and $\beta$ satisfying  the statement of Theorem \ref{main_asym}. 
In particular, by  \cite[Theorem 1.2]{Grill88}, the number $I(L_1^{\phi_k}, L_2^{\phi_k})$ of $\lambda>0$ satisfying \eqref{lambda_L_asym} is estimated by 
\begin{equation}\label{I_estim_assym}
 n(PL_1^{\phi_k})-n(P(L_2^{\phi_k})^{-1})\leq I(L_1^{\phi_k}, L_2^{\phi_k}),\end{equation}
 where $P$ is the orthogonal projection on $\ker(L_2^{\phi_k})^\perp=\{v\in L^2(\Gamma):\, (v, \phi_k)_2=0\}$.  As before, we assume that $L^2(\Gamma)$ is endowed with the usual complex scalar product.
 From inequality \eqref{I_estim_assym} and Proposition \ref{Morse_index_assym} we get:
 
  \vspace{0.3cm}
  
$\bullet$\, $\beta<0\quad\Rightarrow\quad\ I(L_1^{\phi_k}, L_2^{\phi_k})\geq k-1-0=k-1. $  
 
 \vspace{0.5cm}
 
$\bullet$\, $\beta>0\quad\Rightarrow\quad
 I(L_1^{\phi_k}, L_2^{\phi_k})\geq 2N-k-3-1-N+1=N-k-3.$ 

Observe that $"-1"$ in the estimates above appears due to the fact that $P$ is the orthogonal projection onto the subspace of co-dimension  $1$. 
    The proof of the orbital instability for $p>2$ repeats the one in the case of the profile $\phi_\beta$. 
  \end{proof}
    \begin{remark}
 It is impossible to apply the methods used to prove  \cite[Theorem 1.3 and Theorem 1.4]{NataMasa2020} to obtain some  strong instability results for   profiles $\phi_{k}(x)$ different from $\phi_\beta(x)$ and $\phi_{\frac{N}{2}}(x)$ since we do not have variational characterization of them.
\end{remark}
We finish this paper by
  \begin{proof}[Proof of Theorem \ref{min_phi_1}]
  By Proposition \ref{ground}, the solution to  \eqref{d} is a critical point of $S_\omega$. On the other hand, by Theorem \ref{phi_k}, the set of critical points is given by $\phi_\beta,\phi_k, k=1,\ldots, N-1$. We call the minimizer $\phi$. Then we get for $v=v_1+iv_2\in H^1(\Gamma)$
  $$\langle S''_\omega(\phi)v, v\rangle_{(H^1)'\times H^1}= \mathfrak{t}_1^{\phi}(v_1)+ \mathfrak{t}_2^{\phi}(v_2),$$
 where $\mathfrak{t}_j^{\phi}, j=1,2,$ is the quadratic form associated with the self-adjoint operator $L^{\phi}_j$.
 Using  the Implicit Function Theorem (see, for instance, the proof of \cite[Lemma 29]{RJJ}),  it can be proven that $\langle S''_\omega(\phi)v, v\rangle_{(H^1)'\times H^1}\geq 0$ on the subspace $I'=\{v\in H^1(\Gamma): I'_\omega(\phi)v=0\}$ of codimension one.
 Thus, we have  $\mathfrak{t}_1^{\phi}\geq 0$ on the subspace of codimension one.
Indeed, for the real-valued function  $w\in I'$ we have  $\mathfrak{t}_1^{\phi}(w)=\langle S''_\omega(\phi)w, w\rangle_{(H^1)'\times H^1}.$
 Therefore, $L^{\phi}_1$  has at most one negative eigenvalue. In fact, it has one negative eigenvalue since $(L^{\phi}_1\phi,\phi)_2=-(p-1)\|\phi\|^{p+1}_{p+1}<0.$ 
By Proposition \ref{Morse_index}-$(i),1)$ and    Proposition \ref{Morse_index_assym}-$(ii)$, the only candidate for the minimizer is $\phi_1$, and $n(L_1^{\phi_1})=1$.  
  \end{proof}
\textbf{Acknowledgement} I would like to thank the anonymous reviewer whose valuable comments and suggestions helped  to improve this manuscript.

%\bibliographystyle{plain}
%\bibliography{bibliografia} 

\end{document}